\documentclass[12pt]{amsart}
\usepackage{amscd}
\usepackage{amssymb}
\usepackage[pagebackref]{hyperref} 
\usepackage{graphicx}
\usepackage{xcolor}

\def\Null{\operatorname{Null}}

\def\supp{\operatorname{supp}}

\newcommand{\ZZ}{\mathbb Z}
\newcommand{\NN}{\mathbb N}

\newcommand{\C}{\mathcal C}
\newcommand{\QQ}{\mathbb Q}

\newtheorem{lemma}{Lemma}[section]
\newtheorem{corollary}[lemma]{Corollary}
\newtheorem{theorem}[lemma]{Theorem}
\newtheorem{proposition}[lemma]{Proposition}
\newtheorem{definition}[lemma]{Definition}

\newtheorem{example}[lemma]{Example}

\newtheorem{notation}[lemma]{Notation}

\advance\headheight1.15pt

\begin{document}

\title[Robust toric ideals]{On robust toric ideals of weighted oriented graphs}
\author[R. Nanduri]{Ramakrishna Nanduri$^*$}
\address{Department of Mathematics, Indian Institute of Technology
Kharagpur, West Bengal, INDIA - 721302.}
\email{nanduri@maths.iitkgp.ac.in}
\author[T. K. Roy]{Tapas Kumar Roy$^{\dag}$}
\address{Department of Mathematics, Indian Institute of Technology
Kharagpur, West Bengal, INDIA - 721302.}
\email{tapasroy147@kgpian.iitkgp.ac.in, tapasroy147@gmail.com}
\thanks{$^\dag$ Supported by PMRF research fellowship, India.}
\thanks{$^*$ Corresponding author}
\thanks{{\bf AMS Classification 2020:} 13F65, 13A70, 05E40, 05C38, 05C22.}
\maketitle

\begin{abstract}
  In this work, we study the equivalence of various robustness properties of toric ideals of weighted oriented graphs. For any weighted oriented graph $D$, if its toric ideal $I_D$ is generalized robust (or weakly robust), then we show that $D$ does not have forbidden subgraphs $D_1,D_2$ of certain structures. We give a significant class of weighted oriented graphs $D$ whose toric ideals $I_D$ have the following equivalence. 
  \begin{enumerate}
    \item[(i)] $I_{D}$ is strongly robust (equivalently, $I_{D}$ is robust);   
     \item[(ii)] $I_{D}$ is generalized robust (equivalently, $I_{D}$ is weakly robust); 
    \item[(iii)] $D$ does not have subgraphs equal to $D_{1}$ and $D_{2}$.  
\end{enumerate}
  \end{abstract}

\section{Introduction} \label{sec1}

 Let $R=K[x_1,\ldots,x_n]$, where $K$ is a field. Let $M=({\bf x^{a_1}, \ldots, {\bf x^{a_m}}})$ be a monomial ideal in $R$. Let $A=[{\bf a_1} \dots ,{\bf a_m}]$ be an $ n\times m$ matrix whose columns are ${\bf a_1},\ldots {\bf a_m}$. Consider the polynomial ring $K[e_1, \ldots ,e_m]$, define a $K$-algebra homomorphism $\phi$ : $K[e_1, \ldots, e_m]\rightarrow K[x_1, \ldots , x_n]$, by $\phi(e_i) =\textbf{x}^{\textbf{a}_{i}}, 1\le i\le m$. Then the kernel of $\phi$ is called the {\it toric} ideal of $M$ (or $A$) and is denoted by $I_{A}$. We know that $I_A$ is a prime ideal generated by binomials ${\bf e^u}-{\bf e^v}$, where $\deg_A({\bf e^u})=\deg_A({\bf e^v})$, and $\deg_A({\bf e^u}) =u_1{\bf a_1}+\cdots+u_m{\bf a_m}$. 
  
   The general theory of toric ideals has a substantial literature on various algebraic invariants namely, free resolution, Betti numbers, projective dimension etc. Toric ideals are studied by various authors because of their importance in a variety of research fields such as Commutative Algebra, Combinatorics, Algebraic Geometry, Integer Programming, Combinatorial Optimization, Algebraic Coding Theory, Algebraic Statistics, Computational Algebra etc. refer \cite{s95, v95, oh99, ds, es, ms05, oh05, ps14, nn19, nn22, nn23}, and \cite[Chapter 8]{v01}. A binomial ${\bf e^u}-{\bf e^v}$ in $I_A$ is called a {\it primitive} binomial if there exists no other binomial ${\bf e^{u^{\prime}}}-{\bf e^{v^{\prime}}}$ in $I_A$ such that ${\bf e^{u^{\prime}}}\vert {\bf e^u}$ and ${\bf e^{v^{\prime}}} \vert {\bf e^v}$. The set of primitive binomials in $I_A$ is called the {\it Graver basis} of $I_A$ and denoted by $Gr_A$. Note that $I_{A}$ is generated by $Gr_{A}$. The {\it universal Gr\"obner basis} of the toric ideal $I_{A}$, which is denoted by $\mathcal{U}_{A}$, is defined as the union of all reduced Gr\"obner bases with respect to all possible admissible term orders in $K[e_{1},\ldots, e_{m}]$. A minimal binomial generating set of $I_{A}$ corresponds to a minimal Markov basis of $A$ \cite[Theorem 3.1]{ds}. The {\it universal Markov basis} of $A$ is defined as the union of all minimal Markov bases of $A$ and we denote it by $\mathcal{M}_{A}$ \cite[Definition 3.1]{hs07}. An irreducible binomial in $I_{A}$ is called a {\it circuit binomial} if it has minimal support with respect to set inclusion. The set of all circuit binomials in $I_{A}$ is denoted by $\mathcal{C}_{A}$. It is known that $\mathcal{C}_{A}\subseteq \mathcal{U}_{A}\subseteq Gr_{A}$ \cite[Proposition 4.11]{s95}. Also, we know that $\mathcal{M}_{A}\subseteq Gr_{A}$. Recall that a binomial $f\in I_{A}$ is called an {\it indispensable binomial} if $f$ or $-f$ belongs to every minimal binomial generating set of $I_{A}$. A toric ideal $I_{A}$ is called {\it strongly robust} if $Gr_{A}$ is a minimal generating set of $I_{A}$. A toric ideal $I_{A}$ is called  {\it robust} if $\mathcal{U}_{A}$ is a minimal generating set of $I_{A}$. A toric ideal $I_{A}$ is called {\it generalized robust} if $\mathcal{U}_{A}$ equals $\mathcal{M}_{A}$ and $I_{A}$ is called {\it weakly robust} if $Gr_{A}$ equals $\mathcal{M}_{A}$. The various robustness properties of toric ideals have been investigated by many authors for more than two decades. For strongly robust toric ideals $I_A$, we have 
   $Gr_A= \mathcal{U}_A = \mathcal{M}_A$, see \cite[Theorem 4.2]{ptv18}. 
Tatakis proved that a robust toric ideal has a unique minimal system of binomial generators \cite[Theorem 5.10]{t16}. It is a hard problem to check whether a toric ideal is strongly robust, robust, or generalized robust, or weakly robust. Nevertheless explicit computation of universal Gr\"obner basis, Graver basis of a toric ideal is a demanding problem in many research areas, refer \cite{ds, chks06, cdss, hs07}. The robustness and generalized robustness properties of the toric ideals of simple graphs are studied in \cite{bbdlmns15} and \cite{t16}. If $G$ is a simple graph and $I_G$ denote its toric ideal, then   
\begin{enumerate}
  \item  $I_{G}$ is robust $\iff$ $I_{G}$ is strongly robust (\mbox{\cite[Theorem 3.2]{bbdlmns15}}); 
  \item $I_{G}$ is generalized robust $\iff$ $I_{G}$ is weakly robust (\mbox{\cite[Theorem 3.4]{t16}}). 
\end{enumerate}

In this paper, our goal is to study various robustness properties of toric ideals of edge ideals of weighted oriented graphs. Studying edge ideals of WOGs and their algebraic and combinatorial properties is an active research area at present, see for example \cite{mpv17, gmsvv18, hlmrv19, prt19, zxwt19, kblo22} and references therein. Let $D=(V(D), E(D),w)$ be a (vertex) weighted oriented graph (shortly, WOG), where $V(D)$ is the vertex set of $D$. An edge $e\in E(D)$ is an ordered pair $e=(x_i,x_j)$ with $x_i,x_j\in V(D)$ and the orientation of $e$ is from the vertex $x_i$ to the vertex $x_j$. Also, $w:V(D) \rightarrow \NN$ is a weighted function that assigns weights to each vertex of $D$. Let $V(D)=\{x_{1},\ldots,x_{n}\}$ and $E(D)=\{e_{1},\ldots,e_{m}\}$ is the edge set of $D$. Let $w_{i}$ denote the weight of the vertex $x_{i}$. The {\it edge ideal} of $D$ is an ideal in $R$ defined as  
$$I(D)=(x_{i}x_{j}^{w_{j}} : \mbox{there is an edge from}\; x_{i}\; \mbox{to}\; x_{j}).$$

The toric ideal of $D$, denoted by $I_{D}$, is the toric ideal of the edge ideal $I(D)$ of $D$. We denote Graver basis, universal Gr\"obner basis, universal Markov basis, the set of circuit binomials of $I_{D}$ by $Gr_{D}, \mathcal{U}_{D}, \mathcal{M}_{D}, \mathcal{C}_{D}$ respectively. In \cite{nr24}, the authors gave a characterization of $\C_D$ for any weighted oriented graph $D$, in terms of minors of the incidence matrix $A(D)$ of $D$. Recently in \cite{nr24 a}, the authors studied the strongly robust property of $I_D$. In this work, we investigate the equivalence of the strongly robustness and the robustness (or weakly robustness) of $I_D$. Suppose $D=\C_{1}\cup_P \cdots \cup_P \C_{m}$, where each $\C_j$ is a weighted oriented chordless cycle and all these cycles that share a common path $P$. Then we prove that $Gr_{D}=\mathcal{U}_{D}$ and this set is a reduced Gr\"obner basis with respect to degree lexicographic order (Theorem \ref{sec3thm4}). As a consequence, we show that  
\begin{enumerate}
  \item $I_{D}$ is robust $\iff$ $I_{D}$ is strongly robust;  
\item $I_{D}$ is generalized robust $\iff$ $I_{D}$ is weakly robust,  
\end{enumerate}
see Corollary \ref{sec3cor1}. A weighted oriented cycle $\C$ is called {\it balanced} if det$(A(\C))=0$, where $A(D)$ denotes the incidence matrix of $D$. If $\C$ is a balanced chordless cycle, then by \cite[Theorem 6.1]{bklo}, we have that $I_{\C}$ is a principal ideal, say $I_{\C}=(f_{{\bf c}})$. If $f_{\bf c}$ is indispensable, then $\C$ is not equal to the graph $D_1$ (see Figure \ref{fig2}), where $D_{1}=\C_1\cup_e \C_2$, for some balanced cycles $\C_1$ and $\C_2$  that share an edge $e$ (Proposition \ref{sec3pro6}). If ${\C}$ has exactly one chord, then we show that $f_{\bf c}$ is indispensable iff $\C$ is not equal to the graph $D_1$ (Theorem \ref{sec3thm15}). For any weighted oriented graph $D$, if $I_{D}$ is generalized robust (or weakly robust), then we prove that $D$ has no subgraphs equal to $D_1$ and $D_{2}$ as in Figure \ref{fig2}, where $D_1$ is as described above and $D_{2}=\C_{1}\cup_e \C^{\prime} \cup_e \C^{\prime \prime}$, with $\C_{1}$ a balanced cycle, and  $\C^{\prime}, \C^{\prime \prime}$ unbalanced cycles, such that $\C_1, \C^{\prime}, \C^{\prime \prime}$ share a common edge $e$ (Theorem \ref{sec3pro4}). 
Furthermore, if $D=\C_{1}\cup_P \cdots \cup_P \C_{m}$, such that at most two of the cycles $\C_1,\dots,\C_m$ (is) are unbalanced, then the following statements are equivalent  
\begin{enumerate}
    \item[(i)] $I_{D}$ is strongly robust (equivalently, $I_{D}$ is robust);  
     \item[(ii)] $I_{D}$ is generalized robust (equivalently, $I_{D}$ is weakly robust); 
    \item[(iii)] $D$ has no subgraphs equal to $D_{1}$ and $D_{2}$, 
\end{enumerate}
see (Theorem \ref{sec3thm22}). 
\vskip 0.2cm
The paper is organized as follows. In section \ref{sec2}, we recall the definitions and basic results require for the rest of the paper. In section \ref{sec3}, we prove our main results.

\section{Preliminaries} \label{sec2}

Let $R=K[x_1,\ldots,x_n]$, where $K$ is a field. Let $M$ be a monomial ideal in $R$ minimally generated by monomials $\{\textbf{x}^{\textbf{a}_1},\textbf{x}^{\textbf{a}_2}, \ldots ,\textbf{x}^{\textbf{a}_m}\}$. Let $S=K[e_1, \ldots ,e_m]$   be a polynomial ring. Then define a $K$-algebra homomorphism $\phi$ : $K[e_1, \ldots, e_m]\rightarrow K[x_1, \ldots , x_n]$, as $\phi(e_i) =\textbf{x}^{\textbf{a}_{i}}$. Then the kernel of $\phi$ is called the {\it toric} ideal of $M$, and we denote by $I_M$. Then by \cite[2.1]{bklo}, the irreducible binomials $\displaystyle \prod_{k=1}^{m} e_k^{p_k} - \prod_{k=1}^{m} e_k^{q_k}$ generate $I_M$, satisfying $\displaystyle \sum_{k=1}^m p_k\textbf{a}_k =\sum_{k=1}^m q_k\textbf{a}_k$, with $p_k$ and $q_k$ cannot both be non-zero simultaneously, for any $k$.  Recall that the {\it support} of a binomial $ f=\displaystyle \prod_{k=1}^{m} e_k^{p_k} - \prod_{k=1}^{m} e_k^{q_k}$, denoted as $\supp(f)$, is the set of variables  $e_k$ which appear in $f$ with non-zero exponent.   
 
 A {\it (vertex) weighted oriented graph} is a triplet $D= (V(D),E(D),{\bf w})$, where  $V(D) = \{x_1, \ldots,x_n\}$ is the vertex set of $D$, and 
 $$E(D)=\{(x_i,x_j): \mbox{there is an edge from the vertex}\; x_i \mbox{ to the vertex}\; x_j\}$$ 
 is called the edge set of $D$ and the weight function ${\bf w}:V(D)\rightarrow \NN$. We simply denote the weight function ${\bf w}$ by the vector ${\bf w}=(w_1,\ldots ,w_n)$, where $w_j=w(x_j)$ for each $j$. If $e=(x_i,x_j)$ is an edge of $D$, then we say $x_j$ is the {\it head} of $e$ and $x_i$ is the {\it tail} of $e$. The {\it edge ideal} of $D$ is defined as the ideal 
 $$I(D)= \left(x_ix_j^{w_j} : (x_i, x_j)\in E(D) \right)$$
 in $R$. Then the {\it toric ideal} of $D$ is defined as the toric ideal of $I(D)$ and we shortly denote as $I_D$. Let $E(D)=\{e_1,\ldots,e_m\}$. Then the {\it toric algebra} of $D$ is defined as the algebra $K[e_{1},\ldots,e_{m}]/I_D$. 
 Recall that the {\it incidence matrix} of $D$ is an $n \times m$ matrix whose $(i,j)^{th}$ entry $a_{i,j}$ is defined by  
\begin{center}
$a_{i,j}=
\begin{cases}
1\; \mbox { if } \;e_j = (x_i, x_l) \in E(D)\; \mbox{ for some }\;  1 \le l \le n  \\
w_i\; \mbox{ if }\;e_j = (x_l, x_i) \in E(D)\; \mbox{ for some }\; 1 \le l \le n \\ 0\; \mbox{ otherwise, } 
\end{cases}  
$
\end{center}
and we denote by $A(D)$. We denote $\Null(A(D))$, the null space of $A(D)$ over $\QQ$.  Recall that a weighted oriented even cycle $\mathcal{C}_m$ on $m$ vertices is said to be {\it balanced} if $\displaystyle \prod_{k=1}^{m} a_{k,k} =a_{1,m} \prod_{k=1}^{m} a_{k+1,k}$, that is, $\mbox{det}(A(\mathcal{C}_m))=0$, where $A(\mathcal{C}_m)=[a_{i,j}]_{m \times m}$. 

 \vskip 0.2cm

\begin{definition}
For a vector $\mathbf{b} = \left(b_1,b_2, \ldots,b_{\ell}\right) \in \ZZ^{\ell}$, with all $b_i$'s are non-zero, define the corresponding binomial in $K[e_1,\ldots,e_{\ell}]$ as
$f_{\mathbf{b}} := f_{\mathbf{b}}^{+} -f_{\mathbf{b}}^{-},$  where 
$$f_{\mathbf{b}}^{+} = \prod\limits_{b_i > 0}e_i^{b_i}, \;\;\mbox{ and }\;\; f_{\mathbf{b}}^{-} = \prod\limits_{b_i < 0}e_i^{b_i}.$$  
\end{definition}

\begin{definition} \label{sec2def1}
  For ${\bf a, b} \in \ZZ^l$, we define ${\bf a} \prec {\bf b}$, if $\text{supp}({\bf a}_{+})\subseteq \text{supp}({\bf b}_{+})$ and $\text{supp}({\bf a}_{-})\subseteq \text{supp}({\bf b}_{-})$. For binomials $f_{\bf a}, f_{\bf b}$, we define $f_{\bf a} \prec f_{\bf b}$. if ${\bf a} \prec {\bf b}$, equivalently, supp$(f_{\bf a}^+)\subseteq supp(f_{\bf b}^+)$ and 
  supp$(f_{\bf a}^-)\subseteq supp(f_{\bf b}^-)$. 
\end{definition}

For any vector ${\bf b} \in \ZZ^{l}$, we can always write ${\bf b}={\bf b}_{+}-{\bf b}_{-}$ uniquely, for some ${\bf b}_{+}, {\bf b}_{-} \in \ZZ_+^{l}$. Let $[{\bf b}]_{i}$ denotes the $i$-th entry of the vector ${\bf b}$. Define $\mbox{supp}({\bf b})=\{i:[{\bf b}]_{i}\neq 0\}$. For $S\subseteq \mbox{supp}({\bf b})$, ${\bf b}\vert_{S}$ denotes that $\mbox{supp}({\bf b}\vert_{S})=S$ and $[{\bf b}\vert_{S}]_i=[{\bf b}]_i$ for all $i \in S$. Also, a pure binomial $f_{\bf m}\in I_D$ implies that ${\bf m}\in \Null(A(D))$ and $[{\bf m}]_{i}$ denotes the $i$-th entry of ${\bf m}$ corresponding to the edge $e_i\in E(D)$. Recall that the {\it $k$-minors} of an $m\times n$ matrix $A$ are the determinants of submatrices of $A$ of size $k \times k$.
\begin{definition} 
Let $A$ be any $m\times n$ matrix. For any $1\leq k \leq min\{m,n\}$, we denote $M_k(A[i_1,\ldots,i_{m-k}|j_1,\ldots,j_{n-k}])$, the $k^{th}$ minor of $A$ by deleting the rows $i_1,\ldots,i_{m-k}$ and deleting columns $j_1,\ldots,j_{n-k}$ from $A$.   
\end{definition}
Below proposition is graph version of \cite[Proposition 4.13]{s95}. Note that $A(H)$ is a matrix obtained from $A(D)$ by deleting some columns.
\begin{proposition} \label{sec2pro1}
   Let $H$ be a oriented subgraph of a weighted oriented graph $D$ such that $V(D)=V(H)$. Then $I_H=I_D\cap K[e_i: e_i \in E(H)]$.
\end{proposition}
Below we provide a lemma which we use in many proofs.
\begin{lemma}(\cite[Lemma 2.7]{nr24}) \label{sec2lem1}
 Let $D$ be any weighted oriented graph and $f_{\bf m} \neq 0 \in I_D$. Let $v \in V(D)$ of degree $n$. If $(n-1)$ edges of $D$ incident with $v$ are not in $\supp(f_{\bf m})$, then the other edge incident with $v$ is not in $\supp(f_{\bf m})$. Moreover, if the edge $e_i$ incident with $v$ belongs to $\supp(f_{\bf m}^{+})(resp. \supp(f_{\bf m}^{-}))$, then there exists an edge $e_j$ incident with $v$ belongs to $\supp(f_{\bf m}^{-})(resp. \supp(f_{\bf m}^{+}))$.  
\end{lemma}

Let ${\C}$ be a cycle in a WOG $D$. A {\it chord} of a cycle ${\C}$ is an edge $e\in E(D)$ such that $e$ not belongs to the edge set of ${\C}$ but vertices incident with $e$ belong to the vertex set of ${\C}$.   
The length of a path $P$ is defined as the number of edges in $P$ and we denote by $l(P)$. For any WOG $D$ and $H$ any subgraph of $D$, we denote $D\setminus H$ is the subgraph of $D$ having edge set $E(D) \setminus E(H)$.
Also, we know that for any matrix $A$ with entries of non-negative integers, the corresponding toric ideal $I_A$ is generated by the set of primitive binomials belonging to $I_A$. That is, $I_A$ is generated by $Gr_A$.

\begin{theorem} \cite[Theorem 4.4]{nr24} \label{sec4thm1}
Let $\C_{n}$ be a balanced cycle with vertex set $V(\C_{n})=\{v_1,\ldots,v_{n}\}$ and edge set $E(\C_{n})=\{e_1,\ldots,e_{n}\}$, where the edge $e_i$ is incident with the vertices $v_i$ and $v_{i+1}$ with the convention $v_{n+1}=v_1$. Then $I_{\C_{n}}$ is generated by the primitive binomial $f_{\bf{c_{n}}}$,  where 
$${\bf c_{n}}=\frac{1}{d}\Bigg ((-1)^{i+1}M_{n-1}(A(\C_{n})[1|i]) \bigg)_{i=1}^{n} \in \ZZ^{n},$$ 
$M_{n-1}(A(\C_{n})[1|i])$ denotes the $(n-1)$-minor of $A(\C_n)$ obtained by deleting the first row and the $i^{th}$ column of $A(\C_n)$ and $d=\gcd(M_{n-1}(A(\C_{n})[1|i]))_{i=1}^{n}$. 
\end{theorem}

Throughout this paper, we use the following notation. 

\begin{notation} \label{nota1}
Let $D$ be a weighted oriented graph consisting of two cycles $\C_m, \C_n$ as in Figure \ref{fig1} such that these two cycles share a path $P$ of length $k$, where $1\leq k < \min\{m,n\}$. We denote by $D=\C_m\cup_P \C_n$. Let $V(\mathcal{C}_{m})=\{v_1, v_2, \ldots,v_{k+1}, \ldots  v_{m}\}$,  $V(\mathcal{C}_{n})=\{v_1, v_2, \ldots,v_{k+1},v_{m+1} \ldots,  v_{m+n-k-1}\}$ and  $V(P)=\{v_1,v_{2}, \ldots, v_{k+1}\}$. Let ${\C}$ be the cycle of $D$, whose edge set is given by $E({\C})=(E({\C_m})\cup E({\C_n}))\setminus E(P)$ and $A({\C})$ be the incidence matrix of ${\C}$ with respect to usual labelling where $v_{1}$ corresponds to first row and $v_{m+n-k-1}$ corresponds to the last row of $A({\C})$. Note that $A({\C_m}), A({\C_n})$ are submatrices of $A(D)$ with respect to the induced labelling from $D$.  
\vskip 0.2cm
\begin{figure}[h!] \centering \includegraphics[scale=0.40]{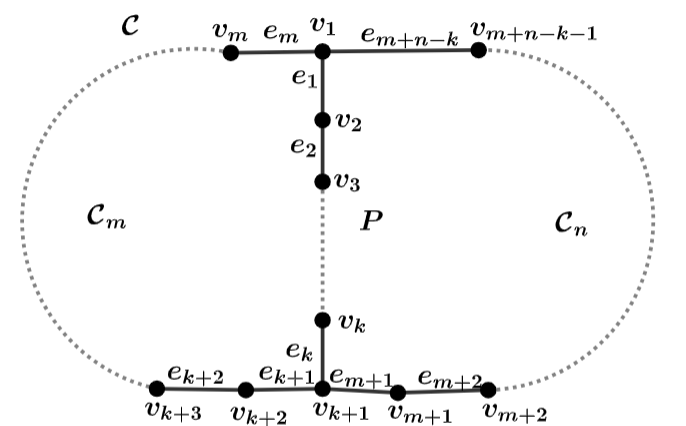}
 \caption{Two weighted oriented cycles that share a path (without orientation)}
    \label{fig1}
\end{figure}
\end{notation}

\begin{notation} \label{sec3nota2}
Let $D$ be a weighted oriented graph. 
\begin{enumerate}
    \item For a balanced cycle $\C_i$ in $D$, we know that its toric ideal $I_{\C_i}$ is generated by a single primitive binomial by \cite[Theorem 4.3]{bklo} and let $f_{\bf c_i}$, denote the primitive generator of $I_{\C_i}$, where ${\bf c_i} \in \Null(A(D))$ is determined by Theorem \ref{sec4thm1}.   
     \item Let ${\C_i}$ and ${\C_j}$ be unbalanced cycles share a path $P$ in $D$. Let $(\C_i \cup \C_j) \setminus P$ denotes the induced subgraph of $D$ whose edge set is $E({\C_i\cup\C_j})\setminus E(P)$. That is, $(\C_i \cup \C_j) \setminus P$ is a cycle and we call it as the outer cycle of $\C_i\cup \C_j$. Assume that the outer cycle of $\C_i\cup \C_j$ is unbalanced. Then by \cite[Theorem 5.1]{bklo}, the toric ideal of $\C_i \cup \C_j$ is principal, say generated by the primitive binomial $f_{\bf{c_ip c_j}}$, for some vector $\bf{c_ip c_j}$ in $\Null(A(D))$. 
\end{enumerate}   
\end{notation}

 

\section{Robustness of toric ideal of weighted oriented graphs} \label{sec3}
In this section, we prove that if $D$ is a weighted oriented graph consisting of certain number of cycles sharing a path, then $Gr_{D}=\mathcal{U}_{D}$ which is a reduced Gr\"obner basis with respect to degree  lexicographic order. Next, we characterize indispensability property of toric ideal of balanced cycle. We show that for any weighted oriented graph $D$, if $I_{D}$ is generalized robust or weakly robust, then $D$ has no subgraphs of types $D_{1}, D_{2}$ as in Figure \ref{fig2}. Also, we prove that if $D$ is a WOG consisting of certain number of cycles sharing a path such that at most two of them are unbalanced and $D$ has no subgraphs of types $D_{1}, D_{2}$, then $I_{D}$ is strongly robust. 

\begin{lemma}\label{lem1}
Let $D=\C_m\cup_P \C_n$ be a weighted oriented graph consisting of two balanced cycles ${\C_m}, {\C_n}$ such that these two cycles share a path $P$ of length $k \ge1$, as in Figure \ref{fig1}. Then the outer cycle ${\C}$ is balanced.
\end{lemma}
\begin{proof}
 Let $A(D)=[b_{i,j}]$ be the incidence matrix of $D$,  where rows are ordered by $v_1, v_2, \ldots, v_{m}, v_{m+2}, \ldots, v_{m+n-k}$ and columns are ordered by $e_1, e_2, \ldots, e_{m+n-k}$. Let $A(\C_{m})$ be the incidence matrix of $\C_{m}$, where rows are ordered by $v_1, v_2, \dots, v_{m}$ and columns are ordered by $e_1, e_2, \dots, e_{m}$. Let $A(\C_{n})$ be the incidence matrix of $\C_{n}$, where rows are ordered by $v_1, v_2, \ldots,v_{k+1}, v_{m+2}, \ldots, v_{m+n-k}$ and columns are ordered by $e_1, e_2, \ldots, e_{k}, e_{m+1}, \ldots, e_{m+n-k}$. Also, assume that $A(\C)$ is the incidence matrix of $\C$, where rows are ordered by $v_{m+1}, v_{m}, \ldots,\\ v_{k+1}, v_{m+2}, \ldots, v_{m+n-k}$ and columns are ordered by $e_{m}, e_{m-1}, \ldots, e_{k+1}, e_{m+1}, \ldots, e_{m+n-k}$. Since $\C_{m}, \C_{n}$ are balanced, we have that det$(A(\C_m))=0$ and det$(A(\C_n))=0$, which implies that 
$$\frac{\prod\limits_{i=1}^{k}b_{i,i}}{\prod\limits_{i=1}^{k}b_{i+1,i}} = \frac{b_{1,m}\prod\limits_{i=k+1}^{m-1}b_{i+1,i}}{\prod\limits_{i=k+1}^{m}b_{i,i}},\;\mbox{and}\; \; \; 
\frac{\prod\limits_{i=1}^{k}b_{i,i}}{\prod\limits_{i=1}^{k}b_{i+1,i}} = \frac{b_{1,m+n-k}\prod\limits_{i=m+1}^{m+n-k-1}b_{i,i}}{b_{k+1,m+1}\prod\limits_{i=m+1}^{m+n-k-1}b_{i,i+1}}.$$
This implies that \begin{equation} \begin{aligned} \label{eq8}\frac{b_{1,m}\prod\limits_{i=k+1}^{m-1}b_{i+1,i}}{\prod\limits_{i=k+1}^{m}b_{i,i}} =\frac{b_{1,m+n-k}\prod\limits_{i=m+1}^{m+n-k-1}b_{i,i}}{b_{k+1,m+1}\prod\limits_{i=m+1}^{m+n-k-1}b_{i,i+1}}.
\end{aligned}
\end{equation} Note that the main diagonal entries of $A(\C)$ are $b_{1,m}, b_{i+1,i}, b_{k+1,m+1}, b_{j,j+1}$ for $i=k+1,\ldots,m-1$,\;  $j=m+1,\ldots,m+n-k-1,$ and the off diagonal entries of $A(\C)$ are $b_{i,i}, b_{1,m+n-k}, b_{j,j}$ for  $i=k+1,\ldots,m,\; j=m+1,\ldots,m+n-k-1$. Thus the product of main diagonal entries of $A(\C)$ is equal to the product of the off diagonal entries of $A(\C)$. 
Now each column of $A({\C})$ has only two non-zero entries, one is on main diagonal and the other is on just below main diagonal and in the last column, off diagonal entry corresponds to the first row. For such a matrix, we have $\mbox{det}(A({\C}))$ is equal to the difference between the product of main diagonal entries and the product of the off diagonal entries. Then by \eqref{eq8}, we get $\mbox{det}(A(\C))=0$. Hence ${\C}$ is balanced.      
\end{proof}
\begin{lemma} \label{sec3lem2}
For any two vectors ${\bf a}$, ${\bf b}\in{\ZZ}^{l}$, suppose $\supp({\bf a})\cap \supp({\bf b})=\{i\}$, and $\lambda_{1}[{\bf a}]_{i}=\lambda_{2}[{\bf b}]_{i}$ for some $\lambda_1,\lambda_2\in \NN$. Then $f_{\lambda_1{\bf a}-\lambda_2{\bf b}}\in(f_{\bf a}, f_{\bf b})$. 
\end{lemma}
\begin{proof}
Without loss of generality, assume that $i\in\mbox{supp}({\bf a}_{+})$. Then by the assumptions, we get $i\in\mbox{supp}({\bf b}_{+})$ and hence $i\notin\mbox{supp}(\lambda_{1}{\bf a}_{}-\lambda_{2}{\bf b}_{})$. Let $S=\mbox{supp}(\lambda_{1}{\bf a}_{}-\lambda_{2}{\bf b}_{})$. Now   
\begin{eqnarray*}
f_{\lambda_1{\bf a}-\lambda_2{\bf b}} &=& f_{(\lambda_1{\bf a}-\lambda_2{\bf b})_{+}}-f_{(\lambda_1{\bf a}-\lambda_2{\bf b})_{-}}\\
& =& f_{(\lambda_1{{\bf a}_{+}}+\lambda_2{{\bf b}_{-}})\vert_{S\setminus\{i\}}} -f_{(\lambda_1{{\bf a}_{-}}+\lambda_2{{\bf b}_{+}})\vert_{S\setminus\{i\}}} \\
  &=& f_{\lambda_{1}{\bf a}_{+}\vert_{\mbox{supp}({\bf a})\setminus \{i\}}}f_{\lambda_2{\bf b}_{-}} - 
  f_{\lambda_1{\bf a}_{-}}f_{\lambda_2{\bf b}_{+}\vert_{\mbox{supp}({\bf b})\setminus \{i\}}} \\
&=& f_{\lambda_1{{\bf a}_{+}}\vert_{\mbox{supp}({\bf a_{}})\setminus \{i\}}} (f_{\lambda_2{\bf b_{-}}}- f_{\lambda_2{\bf b_{+}}}) +f_{\lambda_1{{\bf a}_{+}}\vert_{\mbox{supp}({\bf a_{}})\setminus \{i\}}} f_{\lambda_2{\bf b_{+}}}- f_{\lambda_1{\bf a_{-}}}f_{\lambda_2{{\bf b}_{+}}\vert_{\mbox{supp}({\bf b_{}})\setminus \{i\}}}  
\end{eqnarray*}
 \begin{eqnarray*} 
 &=& f_{{\lambda_1{{\bf a}_{+}}\vert_{\mbox{supp}({\bf a_{}})\setminus \{i\}}}} (f_{\lambda_2{\bf b_{-}}}- f_{\lambda_2{\bf b_{+}}}) +f_{{\lambda_2{{\bf b}_{+}}\vert_{\mbox{supp}({\bf b_{}})\setminus \{i\}}}}(f_{{\lambda_1{{\bf a}_{+}}\vert_{\mbox{supp}({\bf a_{}})\setminus \{i\}}}} e_{i}^{\lambda_2[{\bf b}_{+}]_{i}}- f_{\lambda_1{\bf a_{-}}}) \\ 
 &=& f_{{\lambda_1{{\bf a}_{+}}\vert_{\mbox{supp}({\bf a_{}})\setminus \{i\}}}} (f_{\lambda_2{\bf b_{-}}}- f_{\lambda_2{\bf b_{+}}}) +f_{{\lambda_2{{\bf b}_{+}}\vert_{\mbox{supp}({\bf b_{}})\setminus \{i\}}}}(f_{{\lambda_1{{\bf a}_{+}}\vert_{\mbox{supp}({\bf a_{}})\setminus \{i\}}}} e_{i}^{\lambda_1[{\bf a}_{+}]_{i}}- f_{\lambda_1{\bf a_{-}}}) \\ 
  &=&f_{{\lambda_1{{\bf a}_{+}}\vert_{\mbox{supp}({\bf a_{}})\setminus \{i\}}}} (f_{\lambda_2{\bf b_{-}}}- f_{\lambda_2{\bf b_{+}}}) +f_{{\lambda_2{{\bf b}_{+}}\vert_{\mbox{supp}({\bf b_{}})\setminus \{i\}}}}(f_{\lambda_1{\bf a_{+}}}- f_{\lambda_1{\bf a_{-}}}) \\
 & & \in(f_{\lambda_{1}{\bf a}}, f_{\lambda_{2}{\bf b}})\subseteq (f_{{\bf a}}, f_{{\bf b}}).
\end{eqnarray*}
Thus $f_{\lambda_1{\bf a}-\lambda_2{\bf b}}\in(f_{\bf a}, f_{\bf b})$. 
\end{proof}
\begin{lemma} \label{sec3lem3} 
Let $D=\C_1\cup_P \cdots \cup_P \C_n$ be a weighted oriented graph consisting of $n$ balanced cycles ${\C_{1}}, {\C_{2}},\cdots, {\C_{n}}$ that share a path $P$ of length greater than or equal to $1$. Then $\dim_{\QQ}\Null(A(D))=n$.
\end{lemma}
\begin{proof} First we label vertices and edges of the path $P$ in the order $V(P)=\{v_1,v_2,\ldots,v_{k+1}\}$, $E(P)=\{e_1,e_2,\ldots,e_{k}\}$ where $e_i$ is incident with $v_i$ and $v_{i+1}$ for $i=1,2,\ldots,k$. Next we label vertices and edges of the path ${\C_1}\setminus P$ in the order $V({\C_1}\setminus P)\setminus\{v_1, v_{k+1}\}=\{v_{k+2},\ldots,v_{m_1}\}$, $E({\C_1}\setminus P)=\{e_{k+1},\ldots,e_{m_1}\}$ where $e_i$ is incident with $v_i$ and $v_{i+1}$ for $i=k+1,\ldots,m$ under convention $v_{m_1+1}=v_1$.  In this way, we label vertices and edges of remaining paths in the sequence ${\C_2}\setminus P,{\C_3}\setminus P, \cdots, {\C_{n-1}}\setminus P$ and finally label vertices and edges of the path ${\C_n}\setminus P$ in the order $V({\C_n}\setminus P)\setminus\{v_{1}, v_{k+1}\}=\{v_{i} : i=\sum\limits_{i=1}^{n-1}m_{i}-(n-2)k+2, \ldots,\sum\limits_{i=1}^{n}m_{i}-(n-1)k\}$, $E({\C_n}\setminus P)=\{e_{i} : i=\sum\limits_{i=1}^{n-1}m_{i}-(n-2)k+1, \ldots,\sum\limits_{i=1}^{n}m_{i}-(n-1)k\}$, where $e_{j}$ is incident with $v_{k+1}$ and $v_{j+1}$ for $j=\sum\limits_{i=1}^{n-1}m_{i}-(n-2)k+1$, $e_{i}$ is incident with $v_{i}$ and $v_{i+1}$ for $i=\sum\limits_{i=1}^{n-1}m_{i}-(n-2)k+2,\ldots,\sum\limits_{i=1}^{n}m_{i}-(n-1)k$ under convention $v_{\sum\limits_{i=1}^{n}m_{i}-(n-1)k+1}=v_1$. We take incidence matrix $A(D)$ of $D$ with respect to this labelling of vertices and edges. Let $\widetilde{A(D)}=[a_{i,j}]$ be the matrix obtained from $A(D)$ by deleting columns with respect to last labelling edge of each cycle ${\C_i}$ for $i=1,2,\ldots,n$, that is, $m_{1}$-th column and so on, finally $(\sum\limits_{i=1}^{n}m_{i}-(n-1)k)$-th column. Then for each column $j>1$, there is row $i$ such that $a_{i,j}\neq 0$ but $a_{i,l}=0$ for $l<j$. Thus each column of $\widetilde{A(D)}$ is not a linear combination of its proceeding columns. This implies that columns of $\widetilde{A(D)}$ are linearly independent, where the number of columns of $\widetilde{A(D)}$ is $|E(D)|-n$. Thus the rank of $\widetilde{A(D)}$ is $|E(D)|-n$. Hence rank$(A(D)) \geq |E(D)|-n$. This implies that $\dim_{\QQ} \Null(A(D)) \leq n$. Let $f_{\bf c_i}$ be a primitive binomial that generates $I_{\C_i}$ for $i=1,2,\ldots, n$. Then by the Proposition \ref{sec2pro1}, we have ${\bf c_1, c_2, \ldots, c_n} \in \Null(A(D))$. Since the elements ${\bf c_1, c_2, \ldots, c_n}$ are linearly independent over ${\QQ}$, we get $\dim_{\QQ} \Null(A(D)) \ge n$. Hence $\dim_{\QQ}\Null(A(D))=n$.   
\end{proof}
\begin{lemma} \label{sec3lem4} 
Let $D= \C_1\cup_P \cdots \cup_P \C_n$ be a weighted oriented graph consisting of one unbalanced cycle ${\C_1}$ and $n-1$ cycles ${\C_2},{\C_3},\ldots,{\C_n}$ such that ${\C_1},{\C_2},\ldots,{\C_n}$ share a path $P$ of length greater than or equal to $1$. Then $\dim_{\QQ}\Null(A(D))=n-1$.
\end{lemma}
\begin{proof}
We label the vertices and edges of $D$ as in the proof of the Lemma \ref{sec3lem3}.
Let $A(D^{})$ be the incidence matrix of $D^{}$ with respect to above vertex and edge ordering. We define $\widetilde{A(D^{})}$ is the matrix obtained from $A(D^{})$ by deleting columns corresponding to  last labelling edge of ${\C_j}$ for $j=2,\ldots,n$.  
Then $\widetilde{A(D^{})}$ is a block triangular matrix of the  form   $\left[
\begin{array}{c|c}
  A & B \\
  \hline
  O & C
\end{array}
\right]$, where $A$ is a matrix of size $|V({\C_1})|\times|E({\C_1})|$ matrix and $C$ is an upper triangular matrix of size $(|V(D^{})|-|V({\C_1})|)\times(|E(D^{})|-(n-1)-|E({\C_1})|)$. Note that $A$ is the incidence matrix of ${\C}_{1}$. Since the toric ideal of unbalanced cycle is zero, then det($A$) $\neq$ 0. On the other hand, $C$ is an upper triangular matrix having each diagonal entry is positive. This implies that det($C$) $\neq$ 0. These two arguments show that det($\widetilde{A(D^{})}$) $\neq$ 0. So, the rank of the matrix $\widetilde{A(D^{})}$ is $|V(D^{})|=|E(D^{})|-(n-1)$. Therefore rank$(A(D^{})) \ge |V(D^{})|=|E(D^{})|-(n-1)$. This implies that $\dim_{\QQ}\Null(A(D^{}))\le n-1$. As the number of rows of the matrix $A(D)$ is $|V(D)|$, where $|V(D)|=|E(D)|-(n-1)$, then rank$(A(D^{}))\le |V(D^{})|=|E(D^{})|-(n-1)$. Thus $\dim_{\QQ}\Null(A(D^{}))\ge n-1$. Hence $\dim_{\QQ}\Null(A(D^{}))= n-1$.
\end{proof}

\begin{lemma} \label{sec3lem5}
 Let $D$ be any weighted oriented graph. Let $f_{\bf m}, f_{\bf n}\in I_{D}\subset K[e:e\in E(D)]$ be primitive binomials, where ${\bf m,n} \in \Null(A(D))$. Suppose $e_{i}, e_{j} \in E(D)$ are incident with a vertex $v\in V(D)$ and the degree of $v$ is $2$ in $D$. If $f_{\bf n}^{+}\vert f_{\bf m}^{+}$ and $e_{i}\in\supp(f_{\bf n}^{-})$, then $e_{i}^{[{\bf n_{-}}]_{i}}\vert e_{i}^{[{\bf m_{-}}]_{i}}$. That is, $[{\bf n_{-}}]_{i}\le[{\bf m_{-}}]_{i}$.    
\end{lemma}
\begin{proof}
Since $e_{i}\in\supp(f_{\bf n}^{-})$, by Lemma \ref{sec2lem1}, we have $e_{j}\in\mbox{supp}(f_{\bf n}^{+})$. Since $f_{\bf n}^{+}\vert f_{\bf m}^{+}$, we have $e_{j}\in\mbox{supp}(f_{\bf m}^{+})$ and by Lemma \ref{lem1}, we have $e_{i}\in\mbox{supp}(f_{\bf m}^{-})$. Since $A(D){\bf n}=0, A(D){\bf m}=0$, then comparing rows  corresponding to the vertex $v$, we get $[{\bf n}_{-}]_{i}x=[{\bf n}_{+}]_{j}y$, $[{\bf m}_{-}]_{i}x=[{\bf m}_{+}]_{j}y$  for some positive integers $x, y$. As $f_{\bf n}^{+}\vert f_{\bf m}^{+}$, then $e_{j}^{[{\bf n}_{+}]_{j}}\vert e_{j}^{[{\bf m}_{+}]_{j}}$ and $[{\bf n}_{+}]_{j}\le [{\bf m}_{+}]_{j}$. Then using above equations, we get $[{\bf n}_{-}]_{i}\le [{\bf m}_{-}]_{i}$ and $e_{i}^{[{\bf n}_{-}]_{i}}\vert e_{i}^{[{\bf m}_{-}]_{i}}$.   
\end{proof}

\begin{proposition} \label{sec3pro1}
Let $D=\C_1\cup_P \cdots \cup_P \C_n$ be a weighted oriented graph consisting of $n$ cycles ${\C_1},{\C_2},\ldots,{\C_n}$ that share a path $P$ such that $l(P)\ge 2$, and $l({\C_i}\setminus P)\ge 2$ for $i=1,2,\ldots, n$. Then $I_{D}$ is strongly robust.
\end{proposition}
\begin{proof}
Let $\mathcal{G}(I(D))$ be the minimal monomial generating set of $I(D)$. For any monomial of the form $v_{i}v_{j}^{w_{j}}$ belonging to $\mathcal{G}(I(D))$, either vertex $v_{i}$ or $v_{j}$ is of degree 2. Thus that vertex of degree 2 belongs to exactly support of two monomials belonging to $\mathcal{G}(I(D))$. Then by \cite[Theorem 3.7]{nr24 a}, we get $I_{D}$ is strongly robust. 
\end{proof}

\begin{theorem} \label{sec3thm4}
Let $D=\C_1\cup_P \cdots \cup_P \C_n$ be a weighted oriented graph consisting of chordless cycles ${\C_1},{\C_2},\ldots,{\C_n}$ that share a path $P$. Then $\mathcal{U}_{D}=Gr_{D}$ and it is a reduced Gr\"obner basis for $I_D$ with respect to the degree lexicographic order. 
\end{theorem}
\begin{proof} Let $D$ be the weighted oriented graph consisting of $n$ cycles say ${\C_1},{\C_2},\ldots,{\C_n}$ that share a path $P$. If $l(P)>1$, and $l({\C_i}\setminus P)>1$, for $i=1,2,\ldots,n$, then by Proposition \ref{sec3pro1}, $I_{D}$ is strongly robust and hence $\mathcal{U}_{D}=Gr_{D}$ which is a  reduced Gr\"obner basis with respect to any monomial order by \cite[Theorem 4.2]{ptv18}. Assume that $l(P)=1$ or $l({\C_i}\setminus P)=1$ for some cycle ${\C_i}$. Then we have either $l(P)=1$ or $l({\C_i}\setminus P)=1$ for a unique $i$ because the underlying simple graph has no multiple edges. Then by reviewing there are cycles ${\C_1^{\prime}}$,${\C_2^{\prime}}$,\ldots,${\C_n^{\prime}}$ in $D$ that sharing a path of length $1$ (say $e_1$): if $l(P)=1$, take ${\C_j^{\prime}}={\C_j}$ for $j=1,2,\ldots,n$ and if $l({\C_i}\setminus P)=1$, then take ${\C_i^{\prime}}={\C_i}$ and ${\C_j^{\prime}}=({\C_i}\setminus P)\cup({\C_j}\setminus P)$ for $j\neq i$. As every ideal has a reduced Gr\"obner basis with respect to a monomial order and since $\mathcal{U}_{D}\subseteq Gr_{D}$, we have $Gr_{D}$ is a Gr\"obner basis with respect to any term order. Now we show that $Gr_{D}$ is a reduced Gr\"obner basis with respect to the degree lexicographic order $>$ with $e_1>\mbox{any variable}$. Note that $Gr_{D}$ consists of monic polynomials. Let $LT_{>}(f)$ denote the leading term of any polynomial $f$. Without loss of generality, assume that $LT_{>}(f)=f^{+}$ for any $f\in Gr_{D}$. Let $f_{\bf m}\in Gr_{D}$. Suppose there exists $f_{\bf n}\in Gr_{D}\setminus\{f_{\bf m}\}$ such that $LT_{>}(f_{\bf n})\vert f_{\bf m}^{+}$. If $e_1\in\mbox{supp}(f_{\bf n})$, then $e_1\in\mbox{supp}(LT_{>}(f_{\bf n}))$ and hence $e_1\notin\mbox{supp}(f_{\bf n}^{-})$. Thus $e_1\notin\mbox{supp}(f_{\bf n}^{-})$ and hence supp$(f_{\bf n}^-) \subseteq E(D)\setminus \{e_1\}$. This implies that for any $e_i\in\mbox{supp}(f_{\bf n}^{-})$, there exists $e_j$ such that $e_i,e_j$ are incident with a vertex of degree $2$. Then using Lemma \ref{sec3lem5}, we can show that $f_{\bf n}^{-}\vert f_{\bf m}^{-}$ which is contradiction. Thus $LT_{>}(f_{\bf n})\nmid f_{\bf m}^{+}$ and similarly we can show that  $LT_{>}(f_{\bf n})\nmid f_{\bf m}^{-}$ for every $f_{\bf n}\in Gr_{D}\setminus\{f_{\bf m}\}$. Thus $Gr_{D}$ is a reduced Gr\"obner basis with respect to degree lexicographic order. Therefore $Gr_{D}\subseteq \mathcal{U}_{D}$ and hence $Gr_{D}= \mathcal{U}_{D}$.
\end{proof}

\begin{corollary} \label{sec3cor1}
Let $D=\C_1\cup_P \cdots \cup_P \C_n$ be a weighted oriented graph consisting of chordless cycles $\C_1\ldots,\C_n$ that share a path $P$. Then
 \begin{enumerate}
    \item[(i)] $I_{D}$ is robust iff $I_{D}$ is strongly robust. 
    \item[(ii)] $I_{D}$ is generalized robust iff $I_{D}$ is weakly robust. 
\end{enumerate}
\end{corollary}
\begin{proof}
By Theorem \ref{sec3thm4}, we have $\mathcal{U}_{D}=Gr_{D}$. Then $Gr_{D}$ is a minimal generating set iff $\mathcal{U}_{D}$ is a minimal generating set. Also, $Gr_{D}=\mathcal{M}_{D}$ iff $\mathcal{U}_{D}=\mathcal{M}_{D}$. This completes the proof.  
\end{proof}
\begin{notation}\label{nota2}
(i) Let $D_{1}$ be a weighted oriented graph consisting of two balanced cycles ${\C_1}, {\C_2}$ that share an edge $e$, where ${\C}$ is the cycle with the edge set $E({\C})=(E({\C_1})\cup E({\C_2}))\setminus\{e\}$. Let $D_2$ be the weighted oriented graph consisting of three cycles ${\C_1}, {\C^{\prime}}, {\C^{\prime\prime}}$ that share an edge $\{e\}$ and ${\C_1}$ is balanced, and ${\C^{\prime}}, {\C^{\prime\prime}}$ are unbalanced. Let ${\C^{\prime\prime\prime}}$ be the cycle with the edge set $E({\C^{\prime\prime\prime}})=(E({\C^{\prime}})\cup E({\C^{\prime\prime}}))\setminus \{e\}$ and ${\C^{\prime\prime\prime}}$ be unbalanced. 

\begin{figure}[h!] \centering \includegraphics[scale=0.45]{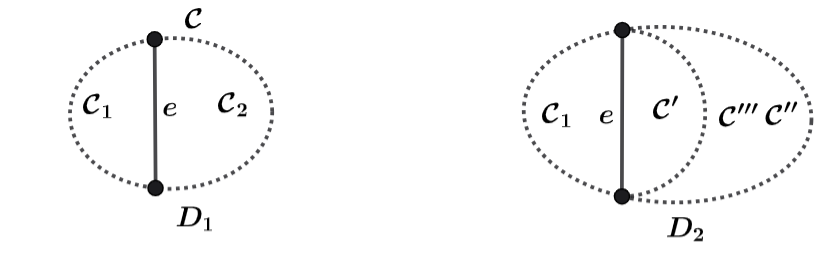}
\caption{Weighted oriented graphs $D_1$ and $D_2$ (without orientation)}
    \label{fig2}
\end{figure} 
\noindent 
(ii) In the graph $D_1$, let $I_{\C_1}=(f_{\bf c_1})$ and $I_{\C_2}=(f_{\bf c_2})$, where $f_{\bf c_1}, f_{\bf c_2}$ are the primitive binomials which can be computed as in Theorem \ref{sec4thm1}. By Lemma \ref{lem1}, the outer cycle $\C$ in $D_1$ is balanced and $I_{\C}=(f_{\bf c})$, where $f_{\bf c}$ can be computed by using Theorem \ref{sec4thm1}. 
\end{notation}

\begin{proposition} \label{sec3pro2}
Let $D_1$ be the graph as in Figure \ref{fig2} and assume Notations \ref{nota2} and \ref{sec3nota2}. Then $f_{\bf c}\in(f_{\bf {c_1}}, f_{\bf {c_2}})$.  
\end{proposition}
\begin{proof}
Assume the Notation \ref{nota2}. We have $I_{\C_1}=(f_{\bf c_1})$ and $I_{\C_2}=(f_{\bf c_2})$. Let $\mathbf{c_1}=\left ((-1)^{i+1}a_i\right )_{i=1}^{m},\mathbf{c_2}=\left (b_1, ((-1)^{i}b_{i})_{i=m+1}^{m+n-1}\right)$, where $a_{i}, b_{i}\in{\NN}$ for each such $i$. Since \\ $\dim_{\QQ} \Null(A(D))=2$ (by Lemma \ref{sec3lem3}), there exist ${\lambda_1^{}},{\lambda_2^{\prime}}\in{\QQ}$ such that ${\bf c}=\lambda_{1}^{}{\bf {c_1}}+\lambda_{2}^{\prime}{\bf {c_2}}$. As $2\in\mbox{supp}({\bf c_1})\cap\mbox{supp}({\bf c}) $, $2\notin\mbox{supp}({\bf c_2})$ and $m+1\in\mbox{supp}({\bf c_2})\cap\mbox{supp}({\bf c}) $, $m+1\notin\mbox{supp}({\bf c_1})$, then the above equation implies that $\lambda_1\neq 0, \lambda_{2}^{\prime}\neq 0$. As $1\notin\mbox{supp}({\bf c})$, we get ${\lambda_1^{}}a_{1}+{\lambda_2^{\prime}}b_1=0$. Then one of $\lambda_{1}, \lambda_{2}^{\prime}$ is positive and the other is negative. Without loss of any generality, assume that $\lambda_{1}>0$, $\lambda_{2}^{\prime}<0$. Let $\lambda_{2}=-\lambda_{2}^{\prime}$ for some $\lambda_{2}>0$. From ${\bf c}=\lambda_{1}^{}{\bf {c_1}}-\lambda_{2}^{}{\bf {c_2}}$, we get $[{\bf c}]_{i}=\lambda_{1}[{\bf c_{1}}]_{i}$ for $i=2,3,\ldots,m$ and $[{\bf c}]_{j}=-\lambda_{2}[{\bf c_{2}}]_{j}$ for $j=m+1,\ldots,m+n-1$. Let $\lambda_{i}=\frac{p_i^{\prime}}{q_i^{\prime}}$,\; $p_i^{\prime},q_i^{\prime}\in{\NN}$, gcd$(p_i^{\prime},q_i^{\prime})=1$ for $i=1,2$. The above equations imply that $q_{1}^{\prime}\vert a_{i}$ for $i=2,\ldots,m$ and $q_{2}^{\prime}\vert b_{j}$ for $j=m+1,\ldots,m+n-1$. Then we can write $\lambda_{i}=\frac{p_i^{\prime\prime}}{d_i^{}}$,\; $p_{i}^{\prime\prime}\in{\NN}$, $i=1,2$ where $d_1=\mbox{gcd}(a_i)_{i=2}^{m}$, $d_2=\mbox{gcd}(b_i)_{i=m+1}^{m+n-1}$. By division algorithm, let $p_{i}^{\prime\prime}=p_{i}+r_{i}d_i,\; 0\le p_{i}<d_{i},\; p_{i},r_{i}\in{\NN}\cup\{0\},\; i=1,2$. Now, we will show that $d_1$ divides $a_1$ and $d_2$ divides $b_1$. 

Let $e$ be the edge incident with vertices $v_{1}, v_{2}$. Suppose $v_1$ is the tail of $e$. Then from $A({\C_1}){\bf {c_1}}=0$, by comparing the row corresponding to the vertex $v_1$, we get $a_1=a_{m}x$ for some positive integer $x$. As $d_1\vert a_{m}$, from above equation, we get $d_1\vert a_{1}$. Suppose $v_2$ is the tail of $e$. Then from $A({\C_1}){\bf {c_1}}=0$, by comparing the row corresponding to the vertex $v_2$, we get $a_1=a_{2}y$ for some positive integer $y$. As $d_1\vert a_{2}$, we have $d_1\vert a_{1}$. Similarly we can show that $d_2\vert b_{1}$. Then $\frac{p_1}{d_1}{\bf {c_1}}\in \Null(A({\C_1}))\cap{\ZZ}^{|E({\C_1})|}$ and hence $f_{\frac{p_1}{d_1}{\bf {c_1}}}\in I_{\C_1}=(f_{\bf {c_1}})$. Suppose $p_1\neq 0$. Then $f_{\frac{p_1}{d_1}{\bf {c_1}}}\neq 0$. As $f_{\frac{p_1}{d_1}{\bf {c_1}}}\in(f_{\bf {c_1}})$, we get ${a}_{i}\le\frac{p_1}{d_1}a_{i}$ for each $i\in\{1,2,\ldots,m\}$ which is contradiction because $\frac{p_1}{d_1}a_{i}<a_{i}$ and  $0<p_{1}<d_{1}$. Thus $p_1=0$. Similarly, $p_2=0$. Then ${\bf c}= -r_{1}{\bf {c_1}}+r_{2}{\bf {c_2}}$, where $r_{1}\in{\NN}, r_{2}\in{\NN}$, $\mbox{supp}({\bf{c_1}})\cap\mbox{supp}({\bf {c_2}})=\{1\}$, $\lambda_{1}[{\bf c_1}]_{1}=\lambda_{2}[{\bf c_2}]_{1}$. Thus using Lemma \ref{sec3lem2}, we get $f_{\bf c}\in(f_{\bf {c_1}}, f_{\bf {c_2}})$.
\end{proof}
\begin{lemma} \label{sec3lem12}
Let ${\C}$ be a balanced cycle in an arbitrary weighted oriented graph $D$ and $I_{\C}=(f_{\bf c})$ for some primitive binomial $f_{\bf c}\in I_{\C}$. Let $f_{\bf n}\in Gr_{D}$ such that $f_{\bf n}^{+}\vert f_{\bf c}^{+}$. If $e_{i}\in\supp(f_{\bf n}^{-})$ and $e_{i}\notin E({\C})$, then $e_{i}$ is a chord of ${\C}$.
\end{lemma}
\begin{proof}
    Suppose $e_{i}$ is not a chord of ${\C}$. Thus there is a vertex, say $v$ incident with $e_{i}$ such that $v\notin V({\C})$. Then by Lemma \ref{sec2lem1}, there is an edge, say $e_{j}\in E(D)$ such that $e_{i}, e_{j}$ are incident at $v$ and $e_{j}\in\mbox{supp}(f_{\bf n}^{+})$. As $v\notin V({\C})$, then $e_{j}\notin E({\C})$ which is a contradiction because $e_{j}\in\mbox{supp}(f_{\bf n}^{+})\subseteq E({\C})$. Thus $e_{i}$ is a chord of ${\C}$.    
\end{proof}
\begin{proposition} \label{sec3lem13}
Let ${\C}$ be a balanced cycle in an arbitrary weighted oriented graph $D$ and $I_{\C}=(f_{\bf c})$ for some primitive binomial $f_{\bf c}\in I_{\C}$. If ${\C}$ has no chord, then $f_{\bf c}$ is indispensable. 
\end{proposition}
\begin{proof}
Suppose $f_{\bf c}$ is not indispensable. Then there is a minimal generating set $M_{D}$ of $I_{D}$ such that $f_{\bf c}\notin M_{D}$ and $-f_{\bf c}\notin M_{D}$. Since $f_{\bf c}\notin M_{D}$, we have $f_{\bf c}$ is a $S$-linear combination of elements of $M_{D}$. Thus there exists $f_{\bf n}\in Gr_{D}\setminus\{ f_{\bf c}\}$ such that $f_{\bf n}^{+}\vert f_{\bf c}^{+}$. We have two possibilities : either $\mbox{supp}(f_{\bf n}^{-})\subseteq E({\C})$ or $\mbox{supp}(f_{\bf n}^{-})\nsubseteq E({\C})$. Suppose $\mbox{supp}(f_{\bf n}^{-})\subseteq E({\C})$. Then $\mbox{supp}(f_{\bf n}^{})\subseteq E({\C})$ and using Proposition \ref{sec2pro1}, $f_{\bf n}\in I_{\C}$. Then either $f_{\bf c}^{+}\vert f_{\bf n}^{+}$ or $f_{\bf c}^{-}\vert f_{\bf n}^{+}$. Without loss of generality, assume that $f_{\bf c}^{+}\vert f_{\bf n}^{+}$. Let $e_{i}\in\mbox{supp}(f_{\bf c}^{+})$. Suppose $e_{i}\in\mbox{supp}(f_{\bf n}^{-})$. As $f_{\bf c}^{+}\vert f_{\bf n}^{+}$, then $e_{i}\in\mbox{supp}(f_{\bf n}^{+})$. Thus $e_{i}\in\mbox{supp}(f_{\bf n}^{+})\cap\mbox{supp}(f_{\bf n}^{-})$ which is contradiction as $\mbox{supp}(f_{\bf n}^{+})\cap\mbox{supp}(f_{\bf n}^{-})=\emptyset$ since $f_{\bf n}$ is pure binomial. Thus $\mbox{supp}(f_{\bf c}^{+})\nsubseteq\mbox{supp}(f_{\bf n}^{-})$. This implies that $f_{\bf c}^{-}\vert f_{\bf n}^{-}$ as $f_{\bf n}$ is multiple of $f_{\bf c}$. We arrive at a contradiction. Thus $\mbox{supp}(f_{\bf n}^{-})\nsubseteq E({\C})$, that is, there is an edge $e_{i}\in\mbox{supp}(f_{\bf n}^{-})$, $e_{i}\notin E({\C})$. Then using Lemma \ref{sec3lem12}, $e_{i}$ is a chord of ${\C}$ which is a contradiction. Hence $f_{\bf c}$ is indispensable.
\end{proof}
\noindent 
The converse of Proposition \ref{sec3lem13} need not be true. See the following example. 

\begin{example}
 Let $D$ be the weighted oriented graph as in Figure \ref{fig3} with the weight vector ${\bf w}=(2,2,3,1,2,2)$. 
\begin{figure}[h!] \centering \includegraphics[scale=0.45]{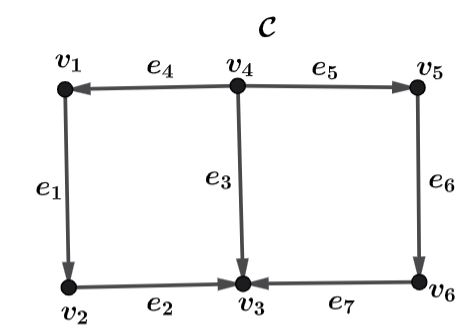}
\caption{Two weighted oriented $4$-cycles that share an edge} \label{fig3}
\end{figure}
\vskip 0.2cm 
Let ${\C}$ be the cycle having edge set $\{e_{1}, e_{2}, e_{7}, e_{6}, e_{5}, e_{4}\}$. Then ${\C}$ is balanced. We see that $I_{D}=(f_{\bf c})$, where $f_{\bf c}=e_{4}e_{2}^{4}e_{6}^{2}-e_{1}^{2}e_{5}e_{7}^{4}$. Thus $f_{\bf c}$ is indispensable but ${\C}$ has chord $e_{3}$.
\end{example}

Using similar arguments like in proof of above Proposition \ref{sec3lem13}, we can prove the following Proposition.

\begin{proposition} \label{sec3lem14}
Let ${\C_1}$, ${\C_2}$ be unbalanced cycles sharing a path $P$ in an arbitrary weighted oriented graph $D$. Let  ${\C_3}$ be the weighted oriented cycle with the edge set $E({\C_3})=(E({\C_1})\cup E({\C_2}))\setminus E(P)$. Suppose ${\C_3}$ is unbalanced. If ${\C_1}, {\C_2}, {\C_3}$ have no chords, then $f_{\bf{c_1pc_2}}$ is indispensable.    
\end{proposition}
\begin{proof}
Since ${\C_1},{\C_2},{\C_3}$ have no chords, then using Lemma \ref{sec2lem1}, it follows that if $f\in Gr_{D}$ such that $f^{+}\vert f_{\bf{c_1pc_2}}^{+}$, then $f=f_{\bf{c_1pc_2}}$. Thus it follows.    
\end{proof}
\begin{proposition} \label{sec3pro5}
Let $D$ be any weighted oriented graph. Suppose $D$ has a subgraph of type $D_{1}$ and assume the notation \ref{nota2}. Then $f_{\bf c}\notin\mathcal{M}_{D}$.   
\end{proposition}
\begin{proof}
 Using Proposition \ref{sec3pro2}, we get $f_{\bf c}\in(f_{\bf c_{1}}, f_{\bf c_{2}})$. Let $M_{D}$ denote any minimal generating set of $I_{D}$. If $f_{\bf c_{1}}\notin M_{D}$, then there is $f_{1}\in M_{D}$ such that $f_{\bf c_{1}}^{+}=f_{1}^{+}h_{1}$ for some monomial $h_{1}$, where $f_{1}\in I_{D^{\prime}}$, $D^{\prime}$ is a weighted oriented subgraph of $D$ consisting of ${\C_{1}}$ and chords of ${\C_{1}}$. To omit $f_{1}^{-}h_{1}$, there is another $f_{2}\in M_{D}$ such that $f_{1}^{-}h_{1}=f_{2}^{+}h_{2}$ for some monomial $h_{2}$, $f_{2}\in I_{D^{\prime}}$. Since ${\C_{1}}$ has finite number of chords, we get that $f_{\bf c_{1}}$ is polynomial combination of elements of $M_{D}\cap I_{D^{\prime}}$. If $f_{\bf c_{2}}\notin M_{D}$, then similarly, we get that $f_{\bf c_{2}}$ is polynomial combination of elements of $M_{D}\cap I_{D^{\prime\prime}}$, where $D^{\prime\prime}$ is a weighted oriented subgraph of $D$ consisting of ${\C_2}$ and chords of ${\C_2}$. But $f_{\bf c}\notin I_{D^{\prime}}+ I_{D^{\prime\prime}}$. Thus $f_{\bf c}\notin M_{D}$ and then $f_{\bf c}\notin \mathcal{M}_{D}$.    
\end{proof}
\begin{proposition} \label{sec3pro6}
 Let ${\C}$ be a balanced cycle in arbitrary weighted oriented graph $D$ and $I_{\C}=(f_{\bf c})$, for some primitive binomial $f_{\bf c}$. If $f_{\bf c}$ is indispensable, then there do not exist two balanced cycles ${\C_1}, {\C_2}$ that share an edge $e$ such that $E({\C})=\{E({\C_1})\cup E({\C_2})\}\setminus \{e\}$.   
\end{proposition}
\begin{proof}
If this is not true, then $D$ has a subgraph of type $D_{1}$. By Proposition \ref{sec3pro5}, $f_{\bf c}\notin\mathcal{M}_{D}$ and then $f_{\bf c}$ is not indispensable which is contradiction. This completes the proof.  
\end{proof}

\noindent 
The converse of Proposition \ref{sec3pro6} need not be true. See the following example. 

\begin{example} \label{eg1}
Let $D$ be the weighted oriented graph as in Figure \ref{fig4} with the weight vector $w=(2,2,3,1,2,2)$. 

\begin{figure}[h!] \centering \includegraphics[scale=0.45]{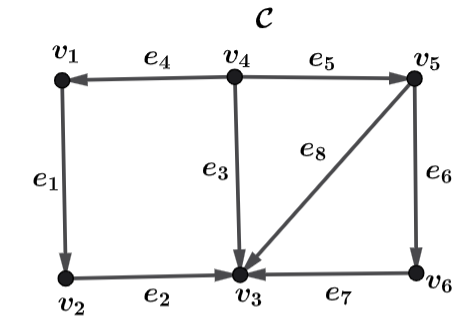}
\caption{}
    \label{fig4}
\end{figure}

Using Macaulay-2 \cite{gs}, We get $Gr_{D}=\{e_{4}e_{2}^{4}e_{6}^{2}-e_{1}^{2}e_{5}e_{7}^{4}, e_{5}e_{7}^{2}-e_{3}e_{6}e_{8}, e_{4}e_{2}^{4}e_{6}-e_{1}^{2}e_{3}e_{7}^{2}e_{8}, e_{4}e_{2}^{4}e_{5}-e_{1}^{2}e_{3}^{2}e_{8}^{2}\}$. Thus $\mathcal{M}_{D}=\{ e_{5}e_{7}^{2}-e_{3}e_{6}e_{8}, e_{4}e_{2}^{4}e_{6}-e_{1}^{2}e_{3}e_{7}^{2}e_{8}, e_{4}e_{2}^{4}e_{5}-e_{1}^{2}e_{3}^{2}e_{8}^{2}\}$. We see that $f_{\bf c}=e_{4}e_{2}^{4}e_{6}^{2}-e_{1}^{2}e_{5}e_{7}^{4}\notin \mathcal{M}_{D}$. Thus $f_{\bf c}$ is not indispensable but there are no balanced cycles ${\C_1}, {\C_2}$ share an edge $e$ such that $E({\C})=\{E({\C_1})\cup E({\C_2})\}\setminus \{e\}$.
\end{example} 

In Example \ref{eg1}, the balanced cycle ${\C}$ has more than one chord. But if ${\C}$ has exactly one chord, then the converse of Proposition \ref{sec3pro6} is true.

\begin{theorem} \label{sec3thm15}
 Let ${\C}$ be a balanced cycle in an arbitrary weighted oriented graph such that ${\C}$ has exactly one chord, say $e$. Let $I_{\C}=(f_{\bf c})$ for some primitive binomial $f_{\bf c}\in I_{\C}$. Then $f_{\bf c}$ is indispensable iff there are no balanced cycles ${\C_1}, {\C_2}$ share an edge $e$ such that $E({\C})=\{E({\C_1})\cup E({\C_2})\}\setminus \{e\}$, that is,  $\C$ is not of the graph type $D_1$.    
\end{theorem}
\begin{proof}
The necessary part follows using Proposition \ref{sec3pro6}.

If possible, let, $f_{\bf c}$ is not indispensable. Then there exists $f_{\bf n}\in Gr_{D}\setminus\{ f_{\bf c}\}$ such that $f_{\bf n}^{+}\vert f_{\bf c}^{+}$. Let $e_{i}\in\mbox{supp}(f_{\bf n}^{-})$ such that $e_{i}\notin E({\C})$. Then by Lemma \ref{sec3lem12}, $e_{i}$ is a chord of ${\C}$. Since ${\C}$ has exactly one chord, we get $e_{i}=e$ and then $\mbox{supp}(f_{\bf n}^{-})\subseteq E({\C})\cup \{e\}$. Since $\mbox{supp}(f_{\bf n}^{+})\subseteq E({\C})$, we get $\mbox{supp}(f_{\bf n}^{})\subseteq E({\C})\cup \{e\}$. Then by Proposition \ref{sec2pro1}, $f_{\bf n}\in I_{D^{\prime}}$, where $D^{\prime}$ consisting of ${\C}$ and chord $e$ of ${\C}$. Let ${\C_1}, {\C_2}$ be the cycles that share edge $e$ such that $E({\C})=\{E({\C_1})\cup E({\C_2})\}\setminus \{e\}$. Then according to the given condition and using Lemma \ref{lem1}, both ${\C_1},{\C_2}$ are unbalanced. Then by \cite[Theorem 5.1, Corollary 5.3]{bklo}, we get $I_{D^{\prime}}=(f_{\bf c})$. This implies that $f_{\bf c}^{+}\vert f_{\bf n}^{+}$ and $f_{\bf c}^{-}\vert f_{\bf n}^{-}$ which is a contradiction. Hence $f_{\bf c}$ is indispensable.   
\end{proof}

\begin{corollary} \label{sec3pro16}
Let $D=\C_1\cup_P \cdots \cup_P \C_n$ be a weighted oriented graph consisting of cycles $\C_1\ldots,\C_n$ that share a path $P$. Let ${\C}$ be a balanced cycle in $D$ and $I_{\C}=(f_{\bf c})$, for some primitive binomial. Then $f_{\bf c}$ is indispensable iff there are no balanced cycles ${\C^{\prime}}, {\C^{\prime\prime}}$ in $D$ sharing an edge $e$ such that $E({\C})=(E({\C^{\prime}})\cup E({\C^{\prime\prime}}))\setminus \{e\}$.   
\end{corollary}
\begin{proof}
The necessary part follows using Proposition \ref{sec3pro6}.

Since the underlying graph has no multiple edges, then ${\C}$ has no chord or ${\C}$ has exactly one chord. If ${\C}$ has no chord, then by Proposition \ref{sec3lem13}, $f_{\bf c}$ is indispensable. Suppose ${\C}$ has exactly one chord. Then by Theorem \ref{sec3thm15}, $f_{\bf c}$ is indispensable.
\end{proof}

\begin{theorem} \label{sec3pro4}
For any weighted oriented graph $D$, if $I_{D}$ is generalized robust or weakly robust, then $D$ has no subgraphs of types $D_1$ and $D_{2}$ as in Figure \ref{fig2}. 
\end{theorem}
\begin{proof}
    Suppose $I_{D}$ is generalized robust or weakly robust. Assume that $D$ has subgraph of type $D_1$ as in Notation \ref{nota2}. Then by Proposition \ref{sec3pro5}, we get $f_{\bf c}\notin\mathcal{M}_{D}$. But $f_{\bf c}\in\mathcal{C}_{D}\subseteq\mathcal{U}_{D}\subseteq Gr_{D}$. Then $\mathcal{U}_{D}\neq\mathcal{M}_{D}$, $Gr_{D}\neq\mathcal{M}_{D}$, that is, $I_{D}$ is not generalized robust and not weakly robust. This is a contradiction. Thus $D$ has no subgraph of type $D_{1}$.

\noindent 
Suppose $D$ has a subgraph of type $D_{2}$ as in Notation \ref{nota2}. Let $P$ be the sharing path of the unbalanced cycles ${\C^{\prime}}, {\C^{\prime\prime}}$. Then $E(P)=\{e_{}\}$.
Let ${\C^{*}}$ be the cycle whose edge set is $(E({\C_1})\cup E({\C^{\prime}}))\setminus E(P)$ and ${\C^{**}}$ be the cycle whose edge set is $(E({\C_1})\cup E({\C^{\prime\prime}}))\setminus E(P)$. Then by Lemma \ref{lem1}, ${\C^{*}}$ and ${\C^{**}}$ are unbalanced. We see that the cycles ${\C^{*}}$, ${\C^{**}}$ share a path, say $P^{\prime}$, where $E(P^{\prime})=E({\C_1})\setminus\{e\}$ such that the outer cycle ${\C^{\prime\prime\prime}}$ is unbalanced. Now we show that $f_{{\bf c^{*}p^{\prime}c^{**}}}\in (f_{{\bf c_1}}, f_{{\bf c^{\prime}pc^{\prime\prime}}})$. Since $\dim_{\QQ}\Null(A(D_2))=2$ (by Lemma \ref{sec3lem4}), there exists $\lambda_1, \lambda_2\in{\QQ}$ such that ${{\bf c^{*}p^{\prime}c^{**}}}=\lambda_1{\bf c_1}+\lambda_2{\bf c^{\prime}pc^{\prime\prime}}$. We can take $({\bf c_1})\vert_{\{i:e_{i}\in E(P^{\prime})\}}\prec({\bf c^{*}p^{\prime}c^{**}})\vert_{\{i:e_{i}\in E(P^{\prime})\}}$ and $({\bf c^{\prime}pc^{\prime\prime}})\vert_{\{i:e_{i}\in E({\C^{\prime}})\setminus E(P)\}}\prec({\bf c^{*}p^{\prime}c^{**}})\vert_{\{i:e_{i}\in E({\C^{\prime}})\setminus E(P)\}}$. Then $\lambda_1>0, \lambda_2>0$. Let $\lambda_i=\frac{p_i^{\prime}}{q_i}$, where $\mbox{gcd}(p_{i}^{\prime}, q_{i})=1$, $p_{i}^{\prime}, q_{i}\in{\NN}$ for $i=1,2$. Then we can write $\lambda_{i}=\frac{p_{i}^{\prime\prime}}{d_i}$, for $i=1,2$, where $d_1=\mbox{gcd}([{\bf c_1}]_{i} : e_i\in E({\C_1})\setminus \{e\})$, $d_2=\mbox{gcd}([{\bf c^{\prime}pc^{\prime\prime}}]_{i} : e_i\in (E({\C^{\prime}})\cup E({\C^{\prime\prime}}))\setminus \{e\})$. Let $p_{i}^{\prime\prime}=p_{i}+r_{i}d_{i}$, for some $0\le p_{i}<d_{i}$, $r_{i}\in{\NN}\cup\{0\}$. Let $e_{}$ be the edge incident with the vertices $v_{1}$ and $v_{2}$ in $D_{2}$. Then either $v_{1}$ or $v_{2}$ is the tail of $e_{}$. Without loss of generality assume that $v_{1}$ is the tail of $e_{}$. Let $e_{i}, e_{j}$ be the other edges of $E({\C^{\prime}})\cup E({\C^{\prime\prime}}) $ such that $e_{i}, e_{j}$ are incident with $v_{1}$. 
Since $A(D_2)({\bf{c^{\prime}pc^{\prime\prime}}})={\bf 0}$, compare the row both sides corresponding to $v_1$, we get $[{\bf c^{\prime}pc^{\prime\prime}}]_{1} +x[{\bf c^{\prime}pc^{\prime\prime}}]_{i} +y[{\bf c^{\prime}pc^{\prime\prime}}]_{j}=0$ for some integers $x, y$. Since $d_{2}\vert[{\bf c^{\prime}pc^{\prime\prime}}]_{i}$, $d_{2}\vert[{\bf c^{\prime}pc^{\prime\prime}}]_{j}$, then the above equation implies that $d_{2}\vert[{\bf c^{\prime}pc^{\prime\prime}}]_{1}$. Thus $\frac{p_2}{d_2}{\bf c^{\prime}p^{}c^{\prime\prime}}\in\mbox{Null}(A(D_2)\cap{\ZZ}^{|E(D_2)|}$ and then $f_{\frac{p_2}{d_2}{\bf c^{\prime}pc^{\prime\prime}}}\in I_{D_2}=(f_{\bf c^{\prime}pc^{\prime\prime}})$. Suppose $p_2\neq 0$. Then $|[{\bf c^{\prime}pc^{\prime\prime}}]_{k}|\le |\frac{p_2}{d_2}[{\bf c^{\prime}pc^{\prime\prime}}]_{k}|$ for each component $k$ which is a contradiction because $|\frac{p_2}{d_2}[{\bf c^{\prime}pc^{\prime\prime}}]_{k}|<|[{\bf c^{\prime}pc^{\prime\prime}}]_{k}|$ as $0<p_1<d_1$. Thus $p_2=0$. Similarly, we can show that $p_1=0$. Then we get ${{\bf c^{*}p^{\prime}c^{**}}}=r_1{\bf c_1}+r_2{\bf c^{\prime}pc^{\prime\prime}}$ for some $r_{1}, r_{2}\in{\NN}$. As $1\notin\mbox{supp}({\bf{c^{*}p^{\prime}c^{**}}})$, we have $r_1[{\bf c_1}]_{1}+r_2[{\bf c^{\prime}pc^{\prime\prime}}]_{1}=0$, that is, $r_1[{\bf c_1}]_{1}-r_2[-({\bf c^{\prime}pc^{\prime\prime}})]_{1}=0$, where $\mbox{supp}({\bf c_1})\cap\mbox{supp}(-{\bf{c^{\prime}pc^{\prime\prime}}})=\{1\}$. Then by Lemma \ref{sec3lem2}, we get $f_{\bf c^{*}p^{\prime}c^{**}}\in(f_{\bf c_1}, f_{\bf c^{\prime}p^{}c^{\prime\prime}})$.  
Let $M_{D}$ denotes any minimal generating set of $I_{D}$. If $f_{\bf c^{\prime}p^{}c^{\prime\prime}}\notin M_{D}$, then using similar arguments like in proof of Proposition \ref{sec3pro5}, we get that $f_{\bf c^{\prime}p^{}c^{\prime\prime}}$ is polynomial combination of elements of $M_{D}\cap I_{D^{\prime}}$ where $D^{\prime}$ is weighted oriented subgraph of $D$ consisting of ${\C^{\prime}},{\C^{\prime\prime}}$, chords of ${\C^{\prime}},{\C^{\prime\prime}},{\C^{\prime\prime\prime}}$. But $f_{\bf c^{*}p^{\prime}c^{**}}\notin I_{D^{\prime}}$. Then we can prove that $f_{\bf{c^{*}p^{\prime}c^{**}}}\notin\mathcal{M}_{D}$. But $f_{\bf{c^{*}p^{\prime}c^{**}}}\in\mathcal{C}_{D}\subseteq\mathcal{U}_{D}\subseteq Gr_{D}$. Thus $\mathcal{U}_{D}\neq \mathcal{M}_{D}, Gr_{D}\neq \mathcal{M}_{D}$, that is, $I_{D}$ is not generalized robust and not weakly robust which is a contradiction. This completes the proof.     
\end{proof}

\noindent 
The converse of Theorem \ref{sec3pro4} need not be true. See Example \ref{sec3example24} below.

\begin{proposition}\label{lem6}
Let $D=\C_1\cup_P \cdots \cup_P \C_n$ be a weighted oriented graph consisting of cycles $\C_1\ldots,\C_n$ that share a path $P$ such that at most one of them is unbalanced. If $D$ has no subgraph of type $D_{1}$, then $I_{D}$ is strongly robust.
\end{proposition}
\begin{proof}
Let ${\C_1},{\C_2},\ldots,{\C_n}$ be balanced cycles in $D$ that share a path $P$. Then according to the given condition, at most one unbalanced cycle share the path $P$ with ${\C_{i}}$ for each $i=1,2,\ldots,n$.
 By the Lemma \ref{lem1}, the cycle in $D$ having edge set $(E({\C_i})\cup E({\C_j}))\setminus E(P)$ is balanced for $i,j\in\{1,2,\ldots,n\}$, $i\neq j$. If $l(P)=1$ or $l({\C_i}\setminus P)=1$ for some $i$, then there are two balanced cycles in $D$ sharing an edge which implies that $D$ has subgraph of type $D_{1}$. Thus $l(P)>1$ and $l({\C_i}\setminus P)>1$ for each $i$. Let $f_{\bf m}\in Gr_{D}$ such that $f_{\bf m}$ not belongs to a minimal generating set of $I_{D}$. Then there exists $f_{\bf t}\in Gr_{D}\setminus\{f_{\bf m}\}$ such that $f_{\bf t}^{+}\vert f_{\bf m}^{+}$. Let $D^{\prime}$ be the graph consisting of balanced cycles ${\C_1},{\C_2},\ldots,{\C_n}$ sharing the path $P$. Since dim$_{\QQ}\Null(A(D))) = n$ (by Lemma \ref{sec3lem4}) and $\{{\bf c_1},{\bf c_2},\ldots,{\bf c_n}\}$ is linearly independent, so $\{{\bf c_1},{\bf c_2},\ldots,{\bf c_n}\}$ is basis of $\Null(A(D))$. Then ${\bf t}$ is a ${\QQ}$-linear combination of ${\bf c_1},{\bf c_2},\ldots,{\bf c_n}$. This implies that $\mbox{supp}(f_{\bf t})\subseteq E(D^{\prime})$. Then for any $e_{i}\in\mbox{supp}(f_{\bf t}^{-})$, there exists an edge $e_{j}\in\mbox{supp}(f_{\bf t}^{+})$ such that $e_{i},e_{j}$ are incident with a vertex of degree 2. Then by Lemma \ref{sec3lem5}, we can show that $f_{\bf t}^{-}\vert f_{\bf m}^{-}$. This yields a contradiction as $f_{\bf m}$ is primitive. Hence $Gr_{D}$ is a minimal generating set of $I_{D}$, that is, $I_{D}$ is strongly robust. 
\end{proof} 

\begin{proposition} \label{lem8}
Let $D= {\C} \cup_P {\C^{\prime}} \cup_P \C_1\cup_P \cdots \cup_P \C_n$ be a weighted oriented graph consisting of cycles ${\C}, {\C^{\prime}}, \C_1\ldots,\C_n$ that share a path $P$ of length at least $1$. Suppose ${\C}, {\C^{\prime}}$ are unbalanced and ${\C_1},{\C_2},\ldots,{\C_n}$ are balanced. If $D$ has no subgraphs of types $D_1$ and $D_2$, as in Notation \ref{nota2}, then $I_{D}$ is strongly robust.  
\end{proposition}
\begin{proof}
Let $f_{\bf m}\in Gr_{D}$ such that $f_{\bf m}$ not belongs to a minimal generating set of $I_{D}$. Then there exists $f_{\bf n}\in Gr_{D}\setminus \{f_{\bf m}\}$ such that $f_{\bf n}^{+}\vert f_{\bf m}^{+}$. We want to show that $f_{\bf n}^{-}\vert f_{\bf m}^{-}$. Let $e_{i}\in\mbox{supp}(f_{\bf n}^{-})$ be arbitrary.

\noindent 
{\em Case 1 :} Suppose there exists an edge $e_{j}$ such that $e_{i}, e_{j}$ are incident at a vertex of degree 2. Then by Lemma \ref{sec3lem5}, the exponent of $e_{i}$ in $f_{\bf n}^{-}$ is less than or equal to the exponent of $e_{i}$ in $f_{\bf m}^{-}$. 

\noindent
{\em Case 2 :} Suppose $E(P)=\{e_{i}\}$. Since $D$ has no subgraph of type $D_{1}$, then $n=1$. Let ${\C^{\prime\prime}}$ be the cycle such that $E({\C^{\prime\prime}})=(E({\C})\cup E({\C^{\prime}}))\setminus E(P)$. Since $D$ has no subgraph of type $D_{2}$, then ${\C}^{\prime\prime}$ must be balanced. Thus $D$ consists of a balanced cycle, say ${\C_{j}}$ for some $j\in\{1,2,\ldots,n\}$ and unbalanced cycles ${\C}, {\C^{\prime}}$ share the path $P$. As $\dim_{\QQ}\Null(A(D))=2$ (by Lemma \ref{sec3lem4}) and ${{\bf c_{j}}, {\bf c^{\prime\prime}}}$ are linearly independent in $\Null(A(D))$, then ${\bf m}$, ${\bf n}$ are ${\QQ}$-linear combinations of ${\bf c_{j}}, {\bf c^{\prime\prime}}$. Also, $\mbox{supp}({\bf c_{j}})\cap \mbox{supp}({\bf c^{\prime\prime}})=\emptyset$. Without loss of  generality, assume that ${\bf c_{j}}\vert_{\{k : e_{k}\in E({\C_{j}}\setminus P)\}}\prec {\bf n}$. As $f_{\bf n}^{+}\vert f_{\bf m}^{+}$, it follows that ${\bf c_{j}}\vert_{\{k : e_{k}\in E({\C_{j}}\setminus P)\}}\prec {\bf m}$. Comparing the  $i^{\mbox{th}}$ components corresponding to $e_{i}$, we get $[{\bf n}]_{i}=\lambda [{\bf c_{j}}]_{i}$, $[{\bf m}]_{i}=\mu [{\bf c_{j}}]_{i}$ for some $\lambda,\mu\in{\QQ}$. Since $e_{i}\in\mbox{supp}(f_{\bf n}^{-})$, then $\lambda\neq 0$. In addition, there exists an edge $e_{k}\in E({\C_{i}}\setminus P)$ such that $e_{k}\in\mbox{supp}(f_{\bf n}^{+})\subseteq \mbox{supp}(f_{\bf m}^{+})$. Then $\mu \neq 0$. Comparing the  $k^{\mbox{th}}$ components corresponding to $e_{k}$, we get $[{\bf n}]_{k}=\lambda [{\bf c_{j}}]_{k}$, $[{\bf m}]_{k}=\mu [{\bf c_{j}}]_{k}$. As ${\bf c_{j}}\vert_{\{k : e_{k}\in E({\C_{j}}\setminus P)\}}\prec {\bf n}$ and ${\bf c_{j}}\vert_{\{k : e_{k}\in E({\C_{j}}\setminus P)\}}\prec {\bf m}$, then $\lambda>0, \mu>0$. Since $f_{\bf n}^{+}\vert f_{\bf m}^{+}$, we get $\lambda\le\mu$. Since $0<\lambda\le\mu$ and $[{\bf n}]_{i}<0$, then using the above equations we get $[{\bf n}]_{i}\ge[{\bf m}]_{i}$. This implies that the exponent of $e_{i}$ in $f_{\bf n}^{-}$ is less than or equal to the exponent of $e_{i}$ in $f_{\bf m}^{-}$.

\noindent
{\em Case 3 :} Let $E({\C_{j}}\setminus P)=\{e_{i}\}$ for some $j\in\{1,2\ldots,n\}$. Any cycle with the edge set $(E({\C_i})\cup E({\C_j}))\setminus E(P)$ is balanced (by Lemma \ref{lem1}) for $i\neq j$. Since $D$ has no subgraph of type $D_{1}$, then $n=1$. Then using similar arguments as in {\em Case 2}, we can show that $e_{i}\in\mbox{supp}(f_{\bf m}^{-})$ and the exponent of $e_{i}$ in $f_{\bf n}^{-}$ is less than or equal to the exponent of $e_{i}$ in $f_{\bf m}^{-}$.

\noindent
{\em Case 4 :} Suppose $E({\C})\setminus E(P)=\{e_{i}\}$. Let $f_{\bf a}$ be a primitive binomial of the toric ideal of the weighted oriented graph comprised with unbalanced cycles ${\C}$, ${\C^{\prime}}$ that share the path $P$. Then by (\cite[Corollary 5.3]{bklo}), $(E({\C^{\prime}})\setminus E(P))\cup \{e_{i}\}\subseteq \mbox{supp}(f_{\bf a})$. Since dim$_{\QQ}\Null(A(D)) = n+1$ (by Lemma \ref{sec3lem4}) and $\{{\bf a},{\bf c_1},{\bf c_2},\ldots,{\bf c_n}\}$ is linearly independent, so $\{{\bf a},{\bf c_1},{\bf c_2},\ldots,{\bf c_n}\}$ is a basis of $\Null(A(D))$. Then we can write ${\bf n}=\lambda{\bf a}+\sum\limits_{j=1}^{n}\lambda_{j}{\bf c_j}, \; {\bf m}=\mu{\bf a}+\sum\limits_{j=1}^{n}\mu_{j}{\bf c_j}$, for some $\lambda,\lambda_{j},\mu,\mu_{j}\in{\QQ}$. Then as $e_{i}\in\mbox{supp}(f_{\bf n}^{-})$, so $\lambda\neq 0$. Also, we get $[{\bf n}]_{j}=\lambda [{\bf a}]_{j}$ for $e_{j}\in E({\C^{\prime}})\setminus E(P)$. Then we have $E({\C^{\prime}})\setminus E(P)\subseteq\mbox{supp}(f_{\bf n})$. As $\mbox{supp}(f_{\bf n}^{+})\subseteq\mbox{supp}(f_{\bf m}^{+})$, then by Lemma \ref{sec2lem1}, we get  $E({\C^{\prime}})\setminus E(P)\subseteq\mbox{supp}(f_{\bf m})$ and ${\bf n}\vert_{\{j:e_{j}\in E({\C^{\prime}})\setminus E(P)\}} \prec {\bf m}\vert_{\{j:e_{j}\in E({\C^{\prime}})\setminus E(P)\}}$. Without loss of generality, assume that ${\bf a}\vert_{\{j:e_{j}\in E({\C^{\prime}})\setminus E(P)\}} \prec {\bf n}\vert_{\{j:e_{j}\in E({\C^{\prime}})\setminus E(P)\}}$. Then ${\bf a}\vert_{\{j:e_{j}\in E({\C^{\prime}})\setminus E(P)\}} \prec {\bf m}\vert_{\{j:e_{j}\in E({\C^{\prime}})\setminus E(P)\}}$. This implies that $\lambda>0, \mu>0$. Also there exists an edge $e_{j}\in E({\C^{\prime}})\setminus E(P)$ such that $j\in\mbox{supp}({\bf n_{+}})\cap\mbox{supp}({\bf m_{+}})$. Then from above equations, corresponding to the edge $e_{j}$, we get $[{\bf n}]_{j}=\lambda[{\bf a}]_{j}$ and $[{\bf m}]_{j}=\mu[{\bf a}]_{j}$. As $f_{\bf n}^{+}\vert f_{\bf m}^{+}$, so $[{\bf n}]_{j}\le[{\bf m}]_{j}$. This implies that $\lambda\le \mu$. Corresponding to the edge $e_i$, we get $[{\bf n}]_{i}=\lambda[{\bf a}]_{i}$ and $[{\bf m}]_{i}=\mu[{\bf a}]_{i}$. Since $i\in\mbox{supp}({\bf n_{-}}), \lambda>0, \mu>0$, we get $i\in\mbox{supp}({\bf m_{-}})$. Since $0<\lambda\le \mu$ and $[{\bf a}]_{i}<0$, we get $[{\bf n}]_{i}\ge[{\bf m}]_{i}$. Thus the exponent of $e_{i}$ in $f_{\bf n}^{-}$ is less than or equal to the  
exponent of $e_{i}$ in $f_{\bf m}^{-}$.

\noindent 
{\em Case 5 :} If $E({\C^{\prime}})\setminus E(P)=\{e_{i}\}$, then using similar arguments as in {\em Case 4}, we can show that the exponent of $e_{i}$ in $f_{\bf n}^{-}$ is less than or equal to the  
exponent of $e_{i}$ in $f_{\bf m}^{-}$.

\noindent
Thus in any case, we get $f_{\bf n}^{-}\vert f_{\bf m}^{-}$ which is a contradiction as $f_{\bf m}\in Gr_{D}$. Hence $Gr_{D}$ is a minimal generating set of $I_{D}$. Thus $I_{D}$ is strongly robust.
\end{proof}

\noindent 
Now we prove the main result of this section. 

\begin{theorem} \label{sec3thm22}
Let $D=\C_1\cup_P \cdots \cup_P \C_n$ be a weighted oriented graph consisting of cycles $\C_1\ldots,\C_n$ that share a path $P$ such that at most two of them are unbalanced. Then the following are equivalent :
\begin{enumerate}
    \item[(i)] $I_{D}$ is strongly robust;  
    \item[(ii)] $I_{D}$ is robust;
    \item[(iii)] $I_{D}$ is generalized robust; 
    \item[(iv)] $I_{D}$ is weakly robust;
    \item[(v)] $D$ has no subgraphs equal to $D_{1}$ and $D_{2}$ as in Figure \ref{fig2}.  
\end{enumerate}
\end{theorem}
\begin{proof}
It is easy to see $(i)\implies (ii)\implies (iii)$.\\
$(iii)\implies (iv) :$
Follows from Corollary \ref{sec3cor1}. 
 \\
$(iv)\implies (v) :$ Follows from Theorem \ref{sec3pro4}.\\
$(v)\implies (i) :$ Follows from Propositions \ref{lem6}, \ref{lem8}.     
\end{proof}

The equivalence in Theorem \ref{sec3thm22} need not be true if we drop the assumption that $D$ has at most two unbalanced cycles. See the following example. 

\begin{example} \label{sec3example24}
 Let $D$ be a weighted oriented graph consisting of three $3$-cycles ${\C_{1}}, {\C_{2}}, {\C_{3}}$ share a path $P$ of length $1$, where $V({\C_{1}})=\{v_{1},v_{2},v_{3}\}, V({\C_{2}})=\{v_{1},v_{2},v_{4}\}, V({\C_{3}})=\{v_{1},v_{2},v_{5}\}, V(P)=\{v_1, v_2\}$ as in Figure \ref{fig5}. Let ${\bf w}=(1,2,3,4,5)$ be the weight vector. 
 
 \begin{figure}[h!] \centering \includegraphics[scale=0.45]{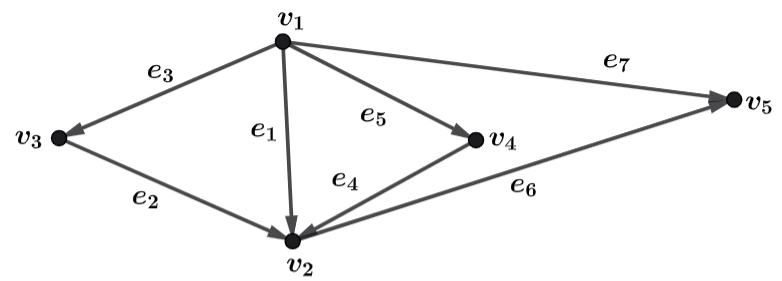}
 \caption{Three weighted oriented $3$-cycles that share an edge}
    \label{fig5}
\end{figure}
We see that every cycle in $D$ is unbalanced. Thus $D$ has no subgraphs of types $D_{1}$ and $D_{2}$. Using Macaulay2 \cite{gs}, we calculate \\ 
$Gr_{D}=\{e_{1}^{3}e_{2}^{3}e_{5}^{2}e_{6}^{4}-e_{3}^{}e_{4}^{8}e_{7}^{4},\;e_{1}^{7}e_{5}^{3}e_{6}^{10}-e_{4}^{12}e_{7}^{10},\;e_{1}^{4}e_{3}^{}e_{5}^{}e_{6}^{6}-e_{2}^{3}e_{4}^{4}e_{7}^{6},\;e_{1}^{1}e_{3}^{2}e_{4}^{4}e_{6}^{2}-e_{2}^{6}e_{5}^{}e_{7}^{2},\;e_{1}^{2}e_{2}^{9}e_{5}^{3}e_{6}^{2}-e_{3}^{3}e_{4}^{12}e_{7}^{2},\;e_{1}^{1}e_{2}^{15}e_{5}^{4}-e_{3}^{5}e_{4}^{16},\;e_{2}^{21}e_{5}^{5}e_{7}^{2}-e_{3}^{7}e_{4}^{20}e_{6}^{2},\;e_{1}^{5}e_{3}^{3}e_{6}^{8}-e_{2}^{9}e_{7}^{8}\}$. \\ 
\noindent 
One can check that every element in $Gr_{D}\setminus\{e_{2}^{21}e_{5}^{5}e_{7}^{2}-e_{3}^{7}e_{4}^{20}e_{6}^{2}\}$ is indispensable binomial of $I_{D}$. Note that  $e_{2}^{21}e_{5}^{5}e_{7}^{2}-e_{3}^{7}e_{4}^{20}e_{6}^{2}=e_{3}^{2}e_{4}^{4}e_{6}^{2}(e_{1}^{}e_{2}^{15}e_{5}^{4}-e_{3}^{5}e_{4}^{16})-e_{2}^{15}e_{5}^{4}(e_{1}^{}e_{3}^{2}e_{4}^{4}e_{6}^{2}-e_{2}^{6}e_{5}^{}e_{7}^{2})$. Thus $Gr_{D}\setminus\{e_{2}^{21}e_{5}^{5}e_{7}^{2}-e_{3}^{7}e_{4}^{20}e_{6}^{2}\}$ is unique minimal generating set of $I_{D}$. Then $\mathcal{M}_{D}=Gr_{D}\setminus\{e_{2}^{21}e_{5}^{5}e_{7}^{2}-e_{3}^{7}e_{4}^{20}e_{6}^{2}\}\neq Gr_{D}$. Hence $I_{D}$ is not weakly robust. Then using Corollary \ref{sec3cor1}, $I_{D}$ is not generalized robust. Thus $I_{D}$ is not robust, not strongly robust.    
\end{example}

The equivalence in Theorem \ref{sec3thm22} can be true without the assumption on $D$. See the following example. 

\begin{example} \label{sec3example23}
 Let $D$ be the weighted oriented graph as in Figure 
\ref{fig6} consisting of three $3$-cycles ${\C_{1}}, {\C_{2}}, {\C_{3}}$ that share an edge $e_1$, where $V({\C_{1}})=\{v_{1},v_{2},v_{3}\}, V({\C_{2}})=\{v_{1},v_{2},v_{4}\}, V({\C_{3}})=\{v_{1},v_{2},v_{5}\}$. Let ${\bf w}=(1,2,3,4,5)$ be the weight vector of vertices of $D$. 
 
 \begin{figure}[h!] \centering \includegraphics[scale=0.45]{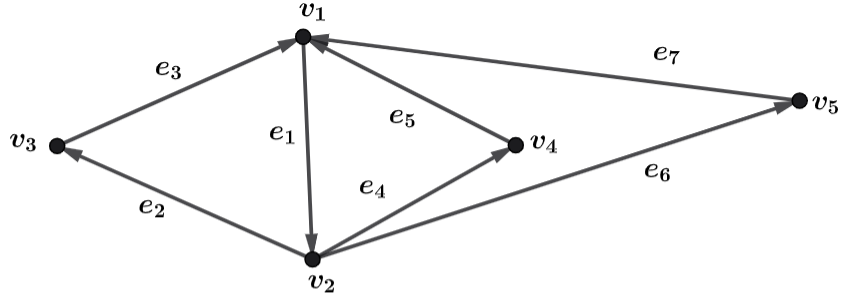}
 \caption{Three weighted oriented $3$-cycles that share an edge}
 \label{fig6}
\end{figure}
Note that every cycle in $D$ is unbalanced. Therefore $D$ has no subgraphs of types $D_{1}$ and $D_{2}$ as in Figure \ref{fig2}. By Macaulay2 \cite{gs}, we have \\ $Gr_{D}=\{e_{2}e_{5}^{8}e_{6}-e_{3}^{3}e_{4}^{2}e_{7}^{5},\; e_{1}e_{5}^{44}e_{6}^{9}-e_{4}^{11}e_{7}^{45},\;e_{1}e_{3}^{3}e_{5}^{36}e_{6}^{8}-e_{2}^{}e_{4}^{9}e_{7}^{40}, \; e_{1}e_{3}^{6}e_{5}^{28}e_{6}^{7}-e_{2}^{2}e_{4}^{7}e_{7}^{35},\\ e_{1}e_{3}^{9}e_{5}^{20}e_{6}^{6}-e_{2}^{3}e_{4}^{5}e_{7}^{30},\;e_{1}e_{3}^{12}e_{5}^{12}e_{6}^{5}-e_{2}^{4}e_{4}^{3}e_{7}^{25},\; e_{1}e_{3}^{15}e_{5}^{4}e_{6}^{4}-e_{2}^{5}e_{4}^{}e_{7}^{20},\; e_{1}e_{3}^{18}e_{4}^{}e_{6}^{3}-e_{2}^{6}e_{5}^{4}e_{7}^{15},\\ e_{1}^{2}e_{3}^{33}e_{6}^{7}-e_{2}^{11}e_{7}^{35},\;e_{1}e_{3}^{21}e_{4}^{3}e_{6}^{2}-e_{2}^{7}e_{5}^{12}e_{7}^{10},\;e_{1}e_{3}^{24}e_{4}^{5}e_{6}^{1}-e_{2}^{8}e_{5}^{20}e_{7}^{5},\;e_{1}e_{3}^{27}e_{4}^{7}-e_{2}^{9}e_{5}^{28}\}$. \\ 
One can check that each binomial in $Gr_{D}$ is indispensable. Therefore $Gr_{D}$ is equal to the set of indispensable binomials of $I_{D}$. Thus $I_{D}$ is strongly robust.
\end{example}

\end{document}